\documentclass[12pt,twoside]{amsart}
\usepackage{amssymb}

\numberwithin{equation}{section}
\newtheorem{theorem}{Theorem}[section]
\newtheorem{lemma}[theorem]{Lemma}
\newtheorem{proposition}[theorem]{Proposition}
\newtheorem{corollary}[theorem]{Corollary}
\newtheorem{conjecture}[theorem]{Conjecture}

\DeclareMathOperator{\gin}{gin}
\theoremstyle{definition}
\newtheorem{definition}[theorem]{Definition} 
\newtheorem{remark}[theorem]{Remark}
\newtheorem{example}[theorem]{Example}

\newtheorem{notation}[theorem]{Notation}

\begin{document}


\newcommand{\m}[1]{\marginpar{\addtolength{\baselineskip}{-3pt}{\footnotesize \it #1}}}
\newcommand{\A}{\mathcal{A}}
\newcommand{\K}{\mathcal{K}} 
\newcommand{\knd}{\mathcal{K}^{[d]}_n}
\newcommand{\F}{\mathcal{F}}
\newcommand{\N}{\mathbb{N}}
\newcommand{\pr}{\mathbb{P}}
\newcommand{\I}{\mathcal{I}}
\newcommand{\G}{\mathcal{G}}
\newcommand{\lcm}{\operatorname{lcm}}
\newcommand{\ndp}{N_{d,p}}
\newcommand{\tor}{\operatorname{Tor}}
\newcommand{\reg}{\operatorname{reg}} 
\newcommand{\mf}{\mathfrak{m}}
\newcommand{\LL}{\mathcal{L}}
\newcommand{\Supp}{\operatorname{Supp}}

 
\title[Componentwise linear monomial ideals]{Some families of componentwise linear monomial ideals}
\thanks{Version: \today}
\author{Christopher A. Francisco}
\address{Department of Mathematics, University of Missouri, 
Mathematical Sciences Building, Columbia, MO 65203}
\email{chrisf@math.missouri.edu}
\urladdr{http://www.math.missouri.edu/$\sim$chrisf}
 
\author{Adam Van Tuyl}
\address{Department of Mathematical Sciences \\
Lakehead University \\
Thunder Bay, ON P7B 5E1, Canada}
\email{avantuyl@sleet.lakeheadu.ca}
\urladdr{http://flash.lakeheadu.ca/$\sim$avantuyl}
 

\begin{abstract}
Let $R=k[x_1,\dots,x_n]$ be a polynomial ring over a field $k$. Let $J=\{j_1,\dots,j_t\}$ be a subset 
of $\{1,\dots,n\}$, and let $\mf_J \subset R$ denote the ideal $(x_{j_1},\dots,x_{j_t})$. Given 
subsets $J_1,\dots,J_s$ of $\{1,\dots,n\}$ and positive integers $a_1,\dots,a_s$, we study ideals of the 
form $I=\mf_{J_1}^{a_1} \cap \cdots \cap \mf_{J_s}^{a_s}$. These ideals arise naturally, for example, 
in the study of fat points, tetrahedral curves, and Alexander duality of squarefree monomial ideals. 
Our main focus is determining when ideals of this form are componentwise linear. Using 
polymatroidality, we prove that $I$ is always componentwise linear when $s \le 3$ or when
$J_i \cup J_j = [n]$ for all $i \neq j$. When 
$s \ge 4$, we give examples to show that $I$ may or may not be
 componentwise linear. We apply these results to ideals 
of small sets of general fat points in multiprojective space, and we extend work of Fatabbi, 
Lorenzini, Valla, and the first author by computing the graded Betti numbers in the $s=2$ case. 
Since componentwise linear ideals satisfy the Multiplicity Conjecture of Herzog, Huneke, and 
Srinivasan when $\operatorname{char}(k) =0$, our work also yields new cases in which this conjecture holds.
\end{abstract}
 
\maketitle

\section{Introduction} \label{s.intro}

Let $R = k[x_1,\ldots,x_n]$ be the polynomial ring in $n$ indeterminates over a field $k$,
and let $[n] :=\{1,\ldots,n\}$.  For a nonempty subset $J = \{j_1,\ldots, j_t\} \subseteq [n]$, we define 
$\mf_J := (x_{j_1},\ldots,x_{j_t})$.  The goal of this paper is to understand when ideals of the form
\[I = \mf_{J_1}^{a_1} \cap \mf_{J_2}^{a_2} \cap \cdots \cap \mf_{J_s}^{a_s}, \, \mbox{with 
$J_i \subseteq [n]$ and $a_i \in \mathbb{Z}^+$},\]
are componentwise linear. We introduce the following definitions.

\begin{definition}\label{d.veronese}
An ideal of the form $\mf_{J_i}^{a_i}$ for some $J_i \subset [n]$ is called a {\bf Veronese ideal} \cite{HHcm}.
We call an ideal $I=\mf_{J_1}^{a_1} \cap \cdots \cap \mf_{J_s}^{a_s}$ an {\bf intersection of Veronese ideals}.
\end{definition}

Let $I \subseteq R$ be a homogeneous ideal, and for a positive integer $d$, let $(I_d)$ be the ideal 
generated by all forms in $I$ of degree $d$. We say that $I$ is {\bf componentwise linear} if for each 
positive integer $d$, $(I_d)$ has a linear resolution. 
 Componentwise 
linear ideals were first introduced by Herzog and Hibi \cite{HH} to generalize Eagon and Reiner's 
result that the Stanley-Reisner ideal $I_{\Delta}$ of a simplicial complex $\Delta$ has a linear 
resolution if and only if the Alexander dual $\Delta^{\star}$ is Cohen-Macaulay \cite{ER}. In 
particular, Herzog and Hibi \cite{HH} and Herzog, Reiner, and Welker \cite{HRW} showed that the 
Stanley-Reisner ideal $I_{\Delta}$ is componentwise linear if and only if $\Delta^{\star}$ is 
sequentially Cohen-Macaulay. On the algebraic side, in characteristic zero, Aramova, Herzog, and Hibi 
subsequently proved that $I$ is componentwise linear if and only if it has the same graded Betti 
numbers as its graded reverse-lex generic initial ideal \cite{AHH}. R\"{o}mer used this result in 
\cite{Roemer} to prove that componentwise linear ideals satisfy the Multiplicity Conjecture of 
Herzog, Huneke, and Srinivasan \cite{HS} in characteristic zero.

Componentwise linearity also arises naturally in the study of several types of ideals from algebraic geometry. 
In \cite{Fra}, the first author showed that if $I$ is the ideal of at most $n+1$ general fat points in 
$\mathbb P^n$, then $I$ is componentwise linear. Additionally, the first author, Migliore, and Nagel proved 
that the ideal of a tetrahedral curve is componentwise linear if and only if the curve does not reduce to a 
complete intersection of type (2,2); see \cite{MN} or \cite{FMN} for an explanation of the reduction process. 
One of our goals in this paper is to identify more results applicable to geometry.

Our motivation to study intersections of Veronese ideals comes from the observation that in many of the cases
 in which the componentwise linear property of a monomial ideal has been studied, the ideal is a special case
of an intersection of Veronese ideals. The defining ideal of $s \leq n$ fat points in $\pr^{n-1}$ in generic 
position, investigated in \cite{Fra}, is an intersection of Veronese ideals with 
$J_i = \{1,\ldots,\hat{i},\ldots,n\}$ for $i=1,\ldots,s$. Moreover, the ideals of tetrahedral curves,
 studied in \cite{MN} and \cite{FMN}, have the form
\[I = (x_1,x_2)^{a_1} \cap (x_1,x_3)^{a_2} \cap (x_1,x_4)^{a_3} \cap 
(x_2,x_3)^{a_4} \cap (x_2,x_4)^{a_5} \cap (x_3,x_4)^{a_6} \subset k[x_1,\dots,x_4],\]
where the $a_i$ are nonnegative integers. Additionally, when each $a_i = 1$, the intersection of Veronese
ideals is the Alexander dual of a Stanley-Reisner ideal; here, the minimal
generators of the Stanley-Reisner ideal are the product of $J_1$ variables, the product of the $J_2$ variables, and so on. 
Faridi showed that if $I$ is the facet ideal
of a simplicial tree (so $I$ is a squarefree monomial
ideal), then the Alexander dual $I^{\star}$ is componentwise linear \cite{Faridi}.
Faridi's result was partially generalized by
the two authors \cite{FVT};  they showed that if $I$ is the edge ideal of a chordal
graph, then the Alexander dual $I^{\star}$ is componentwise linear.

We now present the main results of this paper. Our primary tool is Theorem~\ref{t.maintool}.  We show that if $I$ is an 
intersection of Veronese ideals in $k[x_1,\dots,x_n]$, and if $J_i \cup J_j = [n]$ 
for all $i \not = j$, then $(I_d)$ is a polymatroidal ideal for all $d$. We shall discuss polymatroidal 
ideals in the next section of preliminaries, but their most important property for us is that they have 
linear resolutions. Thus $I$ is componentwise linear in this case since each $(I_d)$
has a linear resolution.  As a corollary of 
Theorem~\ref{t.maintool}, we show that when $s=2$, $I=\mf_J^a \cap \mf_K^b$ is always componentwise linear. With 
some careful analysis of the generators of $(I_d)$, we prove the same result in the case $s=3$ in
 Section~\ref{s.threeideals}.  (When $s=1$, i.e., $I = \mf_J^a$, then the fact that $I$ is componentwise linear 
is simply a corollary of the Eagon-Northcott resolution.)
This shows that the ideals of tetrahedral curves that are not componentwise 
linear given in \cite{FMN} are the simplest possible examples of intersections of Veronese ideals for which 
componentwise linearity fails. When $s \geq 4$, we give examples to show that 
$I = \mf_{J_1}^{a_1} \cap \cdots \cap \mf_{J_s}^{a_s}$ may or may not be componentwise linear. 

In Section~\ref{s.twoideals}, we expand on the $s=2$ case by giving explicit formulas for the graded Betti 
numbers of $\mf_J^a \cap \mf_K^b$. Our formulas generalize results of Fatabbi \cite{Fa}, Valla \cite{Va}, 
Fatabbi and Lorenzini \cite{FL}, and the first author \cite{Fra}, which give the Betti numbers of ideals
of two fat points in $\pr^n$.

We conclude in Section~\ref{s.fatpoints} with some applications. We extend the first author's work in 
\cite{Fra} by showing that if $I$ is the ideal of a small number of general fat points in a multiprojective space 
$\pr^{n_1} \times \cdots \times \pr^{n_r}$, then $I$ is componentwise linear. This also gives a new proof of 
the result in \cite{Fra}; our technique in this paper is more general. Additionally, we use the results of 
Section~\ref{s.twoideals} to write down the graded Betti numbers of two general fat points in multiprojective 
space. We also note that in each case that we show that a class of ideals is componentwise linear, the result 
solves the Multiplicity Conjecture of Herzog, Huneke, and Srinivasan \cite{HS} for that class of ideals 
(in characteristic zero).

\medskip

\noindent
{\bf Acknowledgments.} We gratefully acknowledge the computer algebra systems CoCoA \cite{Co} and Macaulay 2 \cite{M2}, which were invaluable in our work on this paper. The
second author also acknowledges the support provided by NSERC. 
 We also thank Giulio Caviglia for valuable conversations on these topics. Finally, we thank the referee for his or her extremely careful reading of our paper and very helpful corrections and suggestions for improvement.


\section{Preliminaries} \label{s.prelims}

In this section, we recall some definitions and results used throughout the paper. As in the 
introduction, let $R = k[x_1,\ldots,x_n]$ be a polynomial ring over the field $k$,
and for any subset $J = \{j_1,\ldots,j_t\} \subseteq [n]$, we set $\mf_J := (x_{j_1},\ldots,x_{j_t})$.
Our primary interest in this paper is to determine when intersections of Veronese ideals in $R$, 
or equivalently, ideals of the form
$I = \mf_{J_1}^{a_1} \cap \cdots \cap \mf_{J_s}^{a_s},$ where the $a_i$ are 
positive integers,
are componentwise linear.  

Associated to any homogeneous ideal $I$ of $R$ is a {\bf minimal free graded resolution}
\[ 0 \rightarrow \bigoplus_j R(-j)^{\beta_{h,j}(I)} \rightarrow 
\cdots \rightarrow \bigoplus_j R(-j)^{\beta_{1,j}(I)} \rightarrow \bigoplus_j R(-j)^{\beta_{0,j}(I)} 
\rightarrow I \rightarrow 0 \]
where $R(-j)$ denotes the $R$-module obtained by shifting the degrees of $R$ by $j$.  
The number $\beta_{i,j}(I)$ is the $ij$-th graded Betti number of $I$ and equals the number
of generators of degree $j$ in the $i$-th syzygy module.  The following property
of resolutions will be of interest.

\begin{definition} \label{d.linearres}
Suppose $I$ is a homogeneous ideal of $R$ whose generators
all have degree $d$.  Then $I$ has a {\bf linear resolution} if for all $i \ge 0$, $\beta_{i,j}(I) = 0$ 
for all $j \neq i+d$.
\end{definition}

Componentwise linearity is closely related to this property. For a homogeneous ideal $I$, we write $(I_d)$ to 
denote the ideal generated by all degree $d$ elements of $I$.  Note that $(I_d)$ is different from $I_d$, which we 
shall use to denote the vector space of all degree $d$ elements of $I$. Herzog and Hibi introduced the 
following definition in \cite{HH}.
\begin{definition}
\label{d:CWL}
A homogeneous ideal $I$ is {\bf componentwise linear} if $(I_d)$ has a linear 
resolution for all $d$.
\end{definition}

A number of familiar classes of ideals are componentwise linear. For example, all ideals with 
linear resolutions are componentwise linear. However, there are many nontrivial examples as 
well, including stable ideals, squarefree strongly stable ideals, and the $\bold a$-stable ideals 
studied in \cite{GHP}. The following examples illustrates cases in which our results in this paper give new examples of componentwise linear ideals that are not in any of the classes mentioned above.

\begin{example} \label{e.newresults}
Let $R=k[x_1,\dots,x_5]$, and let \[I=(x_1,x_2,x_3) \cap (x_1,x_4,x_5) \cap (x_2,x_3,x_5)=(x_1x_2,x_1x_3,x_1x_5,x_2x_4,x_2x_5,x_3x_4,x_3x_5) \subset R.\] Then $I$ is clearly not stable since no pure power of $x_1$ is among the minimal generators, and it is neither squarefree stable nor $\bold a$-stable for any $\bold a$ because, for example, $x_1x_5$ is a minimal generator, but $x_1x_4 \not \in I$, though other minimal generators do involve $x_4$. Our results in Section~\ref{s.threeideals} show that $I$ is componentwise linear; in fact, $I$ has a linear resolution because it is componentwise linear and has all its minimal generators in the same degree.

For an example that is not squarefree, let $J = (x_1,x_2)^3 \cap (x_2,x_3,x_4,x_5)^2 \subset R$. Then \[ J=({{x}}_{1} {{x}}_{{2}}^{2},{{x}}_{{2}}^{3},{{x}}_{1}^{2} {{x}}_{{2}} {{x}}_{{3}},{{x}}_{1}^{2} {{x}}_{{2}} {{x}}_{{4}},{{x}}_{1}^{2} {{x}}_{{2}} {{x}}_{{5}},{{x}}_{1}^{3} {{x}}_{{3}}^{2},{{x}}_{1}^{3} {{x}}_{{3}} {{x}}_{{4}},{{x}}_{1}^{3} {{x}}_{{3}} {{x}}_{{5}},{{x}}_{1}^{3} {{x}}_{{4}}^{2},{{x}}_{1}^{3} {{x}}_{{4}} {{x}}_{{5}},{{x}}_{1}^{3} {{x}}_{{5}}^{2}),\] which is clearly neither stable nor $\bold a$-stable. By our Theorem~\ref{t.maintool} and Corollary~\ref{c.2CWL}, $J$ is componentwise linear.

\end{example}

The graded Betti numbers of componentwise linear ideals have a particularly good algebraic property. In \cite{AHH}, Aramova, Herzog, and Hibi proved:

\begin{theorem}
\label{t:CWLgin}
Let $I \subset k[x_1,\dots,x_n]$ be a homogeneous ideal, and suppose that $\operatorname{char}(k) = 0$. 
Let $\gin(I)$ be the generic initial ideal of $I$ with respect to the graded reverse-lex order. Then $I$ is componentwise linear if and only if $I$ and $\gin(I)$
have the same graded Betti numbers.
\end{theorem}

In general, $\beta_{i,j}(I) \le \beta_{i,j}(\gin(I))$ for all $i$ and $j$, but all the inequalities are 
equalities exactly when $I$ is componentwise linear. Conca observed in \cite{Conca} that 
Aramova, Herzog, and Hibi actually proved that $I$ is componentwise linear if and only if 
$I$ and $\gin(I)$ have the same number of minimal generators.  This observation makes the condition even easier to test computationally.

One way to show that an ideal is componentwise linear is to prove that it has linear quotients. 
We recall Herzog and Hibi's definition from \cite{HHcm} (which is slightly more restrictive than 
Herzog and Takayama's definition in \cite{HT}).

\begin{definition}
\label{d:linearquo}
Let $I$ be a monomial ideal of $R$. We say that $I$ has {\bf linear quotients} if for some ordering 
$u_1,\dots,u_m$ of the minimal generators of $I$ with 
$\deg u_1 \le \deg u_2 \le \cdots \le \deg u_m$ and all $i > 1$,
 $(u_1,\dots,u_{i-1}):u_i$ is generated by a subset of $\{x_1,\dots,x_n\}$.
\end{definition}

The following proposition is probably known, but we could not find it recorded explicitly, so we
 include it for convenience. The case in which $I$ is generated in a single degree is Lemma 4.1 of \cite{CH}, 
and that is the case we shall use in this paper.

\begin{proposition} \label{p.lq=cwl}
If $I$ is a homogeneous ideal with linear quotients, then $I$ is componentwise linear.
\end{proposition}

\begin{proof}
Suppose that $I \subset R$ has linear quotients with respect to the ordering $u_1,u_2,\dots,u_m$ of its 
minimal generators, where $\deg u_{i-1} \le \deg u_i$ for all $i$. We induct on $m$, the number of minimal generators of $I$. When $m=1$, $I=(u_1)$ is componentwise linear because it is principal.

Fix some $m > 1$. Assume that the ideal $J=(u_1,\dots,u_{m-1})$ is componentwise linear, and suppose that
 $\deg u_m=d$. Let $J'=(J,u_m)$. Note that $J_e=J'_e$ for all $e < d$, so $(J'_e)$ has a linear resolution for 
all $e < d$. We have a short exact sequence
\[ 0 \to R/(J:u_m)(-d) \stackrel{\times u_m}{\longrightarrow} R/J \to R/J' \to 0. \]
Because $J:u_m$ is generated by linear forms,  $\reg(R/(J:u_m))=0$. Since $\deg u_m=d$, we have 
$\reg(R/J') \ge d-1$. Because $R/J$ is componentwise linear, and $\deg u_{m-1} \le d$, we know that 
$\reg(R/J) \le d-1$. By \cite[Corollary 20.19]{E},
\[\reg(R/J') \le \max \{\reg(R/(J:u_n)(-d))-1, \reg(R/J)\}=\max \{d-1,\reg(R/J)\},\] 
so $\reg(R/J')=d-1$. Thus $(J'_d)$ has a linear resolution. The same is true for all $(J'_e)$ with $e > d$. The last statement follows from the fact that for any ideal $M$ with regularity $d$ and $e>d$, $(M_e)$ has a linear resolution.  This fact follows, for example, from \cite[Lemma 2.3]{FMN} since the graded Betti numbers $\beta_{i,j}(M_e)$ with $j > i+e$ must be zero.
\end{proof}

One special type of ideal that has linear quotients is a polymatroidal ideal. For a discussion of this terminology, see \cite{HHdiscrete} and \cite{HHcm}.

\begin{definition}
\label{d:polymatroidal}
Let $I$ be a monomial ideal generated in a single degree. We say that $I$ is a 
{\bf polymatroidal ideal} if the minimal generators of $I$ satisfy the following exchange property: 
If $u=x_1^{a_1} \cdots x_n^{a_n}$ and $v=x_1^{b_1} \cdots x_n^{b_n}$ are minimal generators of $I$, for each $i$ with $a_i > b_i$, there exists $j$ with $a_j < b_j$ such that $x_j u/x_i \in I$.
\end{definition}

Herzog and Takayama proved the following result about polymatroidal ideals in
Lemma 1.3 of \cite{HT}.
                                                                               
\begin{theorem}
\label{t:polymatroidal-lq}
Polymatroidal ideals have linear quotients with respect to the descending
reverse-lex order,
and hence they have linear resolutions.
\end{theorem}
                                                                                
We shall use the ascending reverse-lex order at times, so we state the
corresponding result for that case, which follows from the proof of \cite[Lemma 1.3]{HT} in Herzog and Takayama's paper as well as a dual version of the exchange property for monomial ideals in \cite[Lemma 2.1]{HHcm}.
                                                                                
\begin{proposition}
\label{p.poly-asc}
Polymatroidal ideals have linear quotients with respect to the ascending
reverse-lex order.
\end{proposition}
                                                                                
%
                                                                                Suppose we have a componentwise linear monomial ideal $I=(m_1,\dots,m_r)$ in a polynomial ring 
$R=k[x_1,\dots,x_n]$. In the following sections, we shall sometimes want to consider the ideal 
$I=(m_1,\dots,m_r)$ as an ideal in a larger polynomial ring $R'$. 
The following lemma shows that $I$ is still componentwise linear in the larger ring.

\begin{lemma} \label{l.extendedring}
Let $I=(m_1,\dots,m_r) $ be a componentwise linear monomial ideal in $R=k[x_1,\dots,x_n]$, and let $I'=(m_1,\dots,m_r)R'$ be the ideal generated by the same monomials in the larger polynomial ring $R'=
k[x_1,\dots,x_n,x_{n+1},\dots,x_p]$. Then $I'$ is a componentwise linear ideal of $R'$.
\end{lemma}

\begin{proof} 
Suppose $d$ is the lowest degree in which $I$ has generators. Then $(I_d)=(I'_d)$, so $(I'_d)$ has a linear 
resolution because $(I_d)$ does.

 Let $t \ge 0$, and let $\mf=(x_{n+1},\dots,x_p)$. The ideal $(I'_{d+t})$ has a decomposition as
\[ (I'_{d+t}) = (I_{d+t}) + \mf(I_{d+t-1}) + \mf^2(I_{d+t-2}) + \cdots + \mf^t(I_d);\]
by $(I_{d+u})$, we mean the ideal generated by the degree $(d+u)$ elements of $I$ inside $R$, so the minimal 
generators involve only the variables $x_1,\dots,x_n$. We then consider $\mf^v(I_{d+t-v})$ as an ideal of $R'$.

By hypothesis, $(I_{d+t})$ has a linear resolution in $R$, and hence,
viewed as an ideal of $R'$, we will have $\reg(R'/(I_{d+t})) = d+t-1$. 
We order the rest of the minimal generators of
$(I'_{d+t})$ in the following way. First, take all the minimal generators of $\mf(I_{d+t-1})$ in descending 
graded reverse-lex order (so those monomials divisible by $x_p$ would be last). Next, take all the minimal 
generators of $\mf^2(I_{d+t-2})$ in descending graded reverse-lex order, and continue in this way. 
We shall add each of these 
generators successively to $(I_{d+t})$ and show that each resulting ideal has regularity $d+t$.
This will imply that $\reg(R'/(I'_{d+t})) = d+t-1$ 
and thus $(I'_{d+t})$ has a linear resolution.

As the first step, we compute $(I_{d+t}):x_{n+1}m$, where $m \in I_{d+t-1}$. Multiplying $m$ by any of 
$x_1,\dots,x_n$ gives an element divisible by an element of $I_{d+t}$, and no multiplication by a monomial 
involving only $x_{n+1},\dots,x_p$ can give us an element of $(I_{d+t})$, so $(I_{d+t}):x_{n+1}m=(x_1,\dots,x_n)$. 
We have a short exact sequence
\[0 \rightarrow R'/((I_{d+t}):x_{n+1}m)(-d-t) 
\stackrel{\times x_{n+1}m}{\longrightarrow} R'/(I_{d+t}) \rightarrow R'/((I_{d+t}),x_{n+1}m) \rightarrow 0.\]
By \cite[Corollary 20.19]{E} and the fact that $\reg(R'/(x_1,\ldots,x_n))=0$, we have
\begin{eqnarray*}
\reg(R'/((I_{d+t}),x_{n+1}m)) &\le& \max \{\reg(R'/((I_{d+t}):x_{n+1}m)(-d-t))-1, \reg(R'/(I_{d+t}))\}\\
&=& \max\{d+t-1,d+t-1\} = d+t-1.
\end{eqnarray*}
Since $((I_{d+t}),x_{n+1}m)$ is generated in degree $d+t$, $\reg(R'/((I_{d+t}),x_{n+1}m)) = d+t-1$.

We proceed by induction. Let 
\[ J=(I_{d+t}) + \mf(I_{d+t-1}) + \cdots + \mf^{r-1}(I_{d+t-r+1}) + \, 
\mbox{initial segment of } \, \mf^r(I_{d+t-r}).\]
Suppose $m'=x_{n+1}^{b_{n+1}} \cdots x_p^{b_p} m$ is the next monomial in $\mf^r(I_{d+t-r})$ in descending 
graded reverse-lex order, where $m \in I_{d+t-r}$. First, we will show that $J:m'$ is an ideal 
generated by a 
subset of the variables of $R'$. Multiplying $m$ by any of $x_1,\dots,x_n$ gives an element of $I_{d+t-r+1}$, 
and thus $(x_1,\dots,x_n) \subseteq J:m'$ since 
$x_{n+1}^{b_{n+1}} \cdots x_p^{b_p} \in \mf^r \subset \mf^{r-1}$. Let $l$ be the maximum index for which 
$b_l \not = 0$. Then any of $x_{n+1}m', \dots,x_{l-1}m'$ is in $J$ because 
$x_{n+1} m'/x_l, \dots, x_{l-1}m'/x_l$ are all greater than $m'$ in graded reverse-lex order. 

Now suppose that $\bar{m}$ is a monomial in only $x_l,\dots,x_p$. We will show that 
$\bar{m}m'=\bar{m}x_{n+1}^{b_{n+1}} \cdots x_p^{b_p}m \not \in J$. Note that $\bar{m}m' \in (I_{d+t-r})$, the 
ideal of $R'$ generated by the elements of $I_{d+t-r}$, but it is not in any $(I_{d+t-u})$ for any $u < r$. 
Hence if $\bar{m}m' \in J$, we have $\bar{m}m' \in \mf^r(I_{d+t-r})$. That implies that $\bar{m}m'$ is 
divisible by some monomial in $\mf^r(I_{d+t-r})$ greater than $m'$ in the reverse-lex order. Because of the 
way we have ordered the monomials, and since $l$ is the maximum index for which $b_l \not = 0$, this is 
impossible. Hence 
\[ J:m' = (x_1,\dots,x_n,x_{n+1},\dots,x_{l-1}),\]
an ideal generated by a subset of the variables of $R'$.   We now have an exact
sequence
\[0 \rightarrow R'/(J:m')(-d-t) \stackrel{\times m'}{\longrightarrow} R'/J \rightarrow R'/(J,m') \rightarrow  0.\]
By \cite[Corollary 20.19]{E}, induction, and the fact that $\reg(R'/(J:m'))=0$, we have
\begin{eqnarray*}
\reg(R'/(J,m')) &\le& \max \{\reg(R'/(J:m')(-d-t))-1, \reg(R'/J)\} = d+t-1.
\end{eqnarray*}
Since $(J,m')$ is generated by monomials of degree $d+t$, we have $\reg(R'/(J,m')) = d+t-1$, 
or equivalently, $\reg(J,m') = d+t$ as required.
\end{proof}

\begin{remark} \label{r.withgins}
One can shorten the preceding proof considerably by showing that $\gin(I)$ has the same minimal generators as 
$\gin(I')$, where $\gin$ denotes the graded reverse-lex generic initial ideal. However, this approach would 
require the hypothesis that $\operatorname{char}(k) = 0$ to use the generic initial ideal characterization 
of componentwise linearity.  Instead, we prefer to have a characteristic-free proof.
\end{remark}

We begin our investigation of when intersections of Veronese ideals are componentwise linear with a couple of 
special cases. Let $I = \mf_{J_1}^{a_1} \cap \cdots \cap \mf_{J_s}^{a_s}$. We consider the cases in which 
$s=1$ and in which the $J_i$ are pairwise disjoint.

When $s=1$, $I = \mf_J^a$ is a power of a complete intersection.  In this case, the Eagon-Northcott complex of
 $I$ is a minimal free resolution \cite{EN}. The graded Betti numbers of $I$ are given below (and could also 
be computed from the formulas of \cite{GVT}).

\begin{lemma} \label{betti numbers one J}
Let $J \subseteq [n]$, and let $a$ be any positive integer.  Then
\[\beta_{i,i+a}(\mf_J^a) = \binom{a+|J|-1}{a+i}\binom{a+i-1}{i} ~\mbox{for $i=0,\ldots,|J|-1$},\]
and $\beta_{i,j}(\mf_J^a) = 0$ for all other $i,j \geq 0$.
In particular, $\mf_J^a$ has a linear resolution, and thus is componentwise linear.
\end{lemma}

We use the above lemma to prove the following result.

\begin{theorem}\label{betti numbers many J}
Let $J_1,\ldots,J_s \subseteq [n]$ be $s$ pairwise disjoint nonempty subsets, and let $a_1,\ldots,a_s$ be
positive integers.  Set 
$I = \mf_{J_1}^{a_1} \mf_{J_2}^{a_2} \cdots \mf_{J_s}^{a_s}$ and $a = a_1 + \cdots  + a_s.$  
Then
\[\beta_{i,i+a}(I) = \sum_{i_1+\cdots+i_s = i} \prod_{j=1}^s \binom{a_j + |J_j| -1}{a_j+i_j}\binom{a_j+i_j -1}{i_j} ~~\mbox{for all $i \geq 0$,}\]
and $\beta_{i,j}(I) = 0$ otherwise.
\end{theorem}

\begin{proof}
Let ${\bf F_{\ell}}$ denote the graded minimal free resolution of $\mf_{J_{\ell}}^{a_{\ell}}$  for 
$\ell = 1,\ldots,s$.  
Since $I \cong \mf_{J_1}^{a_1} \otimes \cdots \otimes \mf_{J_s}^{a_s}$, 
the graded minimal free resolution of $I$ is given by
${\bf G} = {\bf F_1} \otimes \cdots \otimes {\bf F_s}$.
So $G_i$, the $i$-th graded free module in a minimal graded free resolution of $I$, is
$G_i = \bigoplus_{i_1+\cdots+i_s =i} F_{i_1}\otimes \cdots \otimes F_{i_s}.$  Thus
\begin{eqnarray*}
\beta_{i,j}(I) &=& \sum_{i_1 + \cdots + i_s = i} \dim_k (F_{i_1} \otimes \cdots \otimes F_{i_s})_j
= \sum_{i_1 + \cdots + i_s = i} \sum_{j_1 + \cdots + j_s = j} 
\beta_{i_1,j_1}(\mf_{J_1}^{a_1}) \cdots \beta_{i_s,j_s}(\mf_{J_s}^{a_s}) 
\end{eqnarray*}
But by Lemma \ref{betti numbers one J}, $\beta_{i_{\ell},j_{\ell}}(\mf_{J_{\ell}}^{a_{\ell}}) \neq 0$
only if $j_{\ell}  = i_{\ell} + a_{\ell}$.  So $\beta_{i,j}(I) = 0$ if $j \neq i+a$, and
\[\beta_{i,i+a}(I) =   \sum_{i_1 + \cdots + i_s = i} 
\beta_{i_1,i_1+a_1}(\mf_{J_1}^{a_1}) \cdots \beta_{i_s,i_s + a_s}(\mf_{J_s}^{a_s}). \]
By applying the  formula of Lemma \ref{betti numbers one J} we get the desired conclusion.
\end{proof}

When the $J_i$'s are pairwise disjoint nonempty sets as in the above theorem,
then  $\mf_{J_1}^{a_1} \cap \mf_{J_2}^{a_2} \cap \cdots \cap \mf_{J_s}^{a_s}
= \mf_{J_1}^{a_1} \mf_{J_2}^{a_2} \cdots \mf_{J_s}^{a_s}$.  Since this
ideal has a linear resolution, we have:

\begin{corollary}If  
$I = \mf_{J_1}^{a_1} \cap \mf_{J_2}^{a_2} \cap \cdots \cap \mf_{J_s}^{a_s}$, with  
the $J_1,\ldots,J_s \subseteq [n]$ pairwise disjoint nonempty subsets, then $I$ is
componentwise linear.
\end{corollary}

\begin{remark}
\label{r:product}
It is easy to see that $(x_1,\dots,x_r)^a$ is polymatroidal for any positive integers $r$ and $a$. Results in \cite{HHdiscrete} and 
\cite{CH} prove that the product of polymatroidal ideals is polymatroidal. Hence the ideals of 
Theorem~\ref{betti numbers many J} are polymatroidal and thus have a linear resolution, 
as is clear from the graded 
Betti numbers.
\end{remark}

\begin{example} We show that if $I = \mf_{J_1}^{a_1} \cap \cdots
\cap \mf_{J_s}^{a_s}$ with $s \geq 4$, then $I$ may or may not be componentwise linear.
First, we construct examples of ideals that are not componentwise linear.  We begin
with the case that $s=4$.
It was observed in \cite{FMN} that the ideal
\[I = (x_1,x_2) \cap (x_2,x_3)
\cap (x_3,x_4) \cap (x_4,x_1) = (x_1x_3,x_2x_4) \]
is not componentwise linear.  To see this fact, note that the ideal 
$I$
is a complete intersection ideal of type $(2,2)$.  Since $I = (I_2)$, $(I_2)$ does not
have a linear resolution.

We can extend this example to  any $s > 4$.  
In the polynomial ring $R = k[x_1,\ldots,x_s]$, let 
\[I =  (x_1,x_2) \cap (x_2,x_3)
\cap (x_3,x_4) \cap (x_4,x_1) \cap (x_5)^{a_5} 
\cap (x_6)^{a_6} \cap \cdots \cap (x_s)^{a_s}. \]
for any positive integers $a_5,\ldots,a_s$.  Then $I = x_5^{a_5}x_6^{a_6}\cdots x_s^{a_s}I'$
where $I'=(x_1x_3,x_2x_4)$.  
Because $\beta_{i,j}(I) = \beta_{i,j-a_5-a_6-\cdots-a_s}(I')$,
the ideal $I$ cannot be componentwise linear since $I = (I_{2+a_5+\cdots+a_s})$ does
not have a linear resolution.

On the other hand, we can create very simple intersections of Veronese ideals that are componentwise linear 
for any $s$. For example, if $J_i = \{i\}$ for $i=1,\dots,s$, then $I$ is principal and hence has a linear 
resolution. Alternatively, start with a componentwise linear intersection of Veronese ideals $I$ in the 
variables $x_1,\dots,x_r$, and intersect $I$ with $(x_{r+1})^{a_{r+1}} \cap \cdots \cap (x_{s})^{a_s}$.
\end{example}

In the following sections, we consider the cases in which $s=2$ or $s=3$ as well as some special cases for general $s$.


\section{A family of polymatroidal ideals} \label{s.poly}

In this section, we consider a particular family of intersections of Veronese ideals. 
We show that ideals in this family are polymatroidal. Our main result is the following theorem.

\begin{theorem} \label{t.maintool}
Let $J_1, \dots, J_s$ be subsets of $[n]$ such that $J_i \cup J_j = [n]$ 
for all $i \not = j$. Let 
\[ I = \mf_{J_1}^{a_1} \cap \cdots \cap \mf_{J_s}^{a_s} \subset R = 
k[x_1,\dots,x_n]. \]
Then $(I_d)$ is polymatroidal for all $d$, and hence $I$ is componentwise 
linear.
\end{theorem}

\begin{proof}
The condition on $J_i \cup J_j = [n]$ means that any $r \in [n]$ is missing 
from at most one of the $J_i$; if $r \not \in J_i$ and $r \not \in J_j$, 
then $J_i \cup J_j \not = [n]$. Therefore we may partition the variables 
$x_1, \dots, x_n$ in the following way: Rename the variables $x_i$ with 
the symbols $x_{1,1},\dots,x_{1,b_1}, \dots,$ $x_{s,1},\dots,x_{s,b_s}$, 
$x_{\cap,1}, \dots, x_{\cap,b_{\cap}}$. The variables $x_{i,j}$ correspond 
to the integers in $[n]$ {\bf missing} from $J_i$, and the variables 
$x_{\cap, j}$ correspond to the integers in $[n]$ present in all $J_i$.

For example, if 
\[ I = (x_1,x_2,x_4,x_6)^3 \cap (x_1,x_3,x_5,x_6)^4 \cap 
(x_2,x_3,x_4,x_5,x_6)^2 \subset k[x_1,\dots,x_6],\]
then $J_1=\{1,2,4,6\}$, $J_2=\{1,3,5,6\}$, and $J_3=\{2,3,4,5,6\}$. We 
would rename the variables $x_3$ and $x_5$ as $x_{1,1}$ and $x_{1,2}$ 
since 3 and 5 are missing from $J_1$. Similarly, $x_2$ and $x_4$ become 
$x_{2,1}$ and $x_{2,2}$, $x_1$ is $x_{3,1}$, and $x_6$ is $x_{\cap,1}$. 
Note that there may be some $i$ with $1 \le i \le s$ for which there are 
no $x_{i,j}$ variables; that is true if and only if $J_i = [n]$. That 
causes no problem in the proof below; alternatively, one can avoid this 
case since a component of $(x_1,\dots,x_n)^a$ simply makes the ideal 
formed by the intersection of the other components zero in degrees below 
$a$ and the same in degrees $a$ and above.

Fix a degree $d$. Suppose that $m_e \not = m_f$ are two monomials in 
$(I_d)$ with 
\[ m_e=x_{1,1}^{e_{1,1}} \cdots x_{1,b_1}^{e_{1,b_1}} \cdots 
x_{s,1}^{e_{s,1}} \cdots x_{s,b_s}^{e_{s,b_s}} x_{\cap,1}^{e_{\cap,1}} 
\cdots x_{\cap,b_{\cap}}^{e_{\cap,b_{\cap}}},\] with $m_f$ having a 
similar expression in terms of $x_{i,j}^{f_{i,j}}$. We need to show that 
the polymatroidal exchange condition holds for these two monomials. 
Namely, if some $e_{i,j} > f_{i,j}$ or some $e_{\cap,j} > f_{\cap,j}$, we 
must show that there exists $e_{u,v} < f_{u,v}$ (with $1 \le u \le s$ or 
$u=\cap$) such that $x_{u,v} m_e / x_{i,j} \in (I_d)$. Note that the fact 
that $m_e \in (I_d)$ means exactly that $\deg m_e = d$ and all of the 
following inequalities hold:
\begin{eqnarray*}
 \sum_j e_{2,j} + \cdots + \sum_j e_{s,j} + \sum_j e_{\cap,j} &\ge &a_1,\\
 \sum_j e_{1,j} + \sum_j e_{3,j} + \cdots + \sum_j e_{s,j} + \sum_j 
e_{\cap,j} &\ge& a_2, \\
\vdots &&\\
 \sum_j e_{1,j} + \cdots + \sum_j e_{s-1,j} + \sum_j e_{\cap,j} &\ge& a_s.
\end{eqnarray*}

There are two main cases to consider. First, suppose that some $e_{\cap,p} 
> f_{\cap,p}$. If there exists $e_{\cap,j} < f_{\cap,j}$, then $x_{\cap,j} 
m_e / x_{\cap,p} \in (I_d)$ since none of the left-hand sides of the 
inequalities above change, and we are done. Otherwise, we have 
$e_{\cap,j} \ge f_{\cap,j}$ for all $j$, and $\sum e_{\cap,j} > 
\sum f_{\cap,j}$ because $e_{\cap,p} > f_{\cap,p}$. Since $m_e$ and $m_f$ have the same degree, there exists 
some $e_{i,j} < f_{i,j}$ with $1 \le i \le s$. Without loss of generality, 
assume that $e_{1,1} < f_{1,1}$. If 
\[ \sum_j e_{2,j} + \cdots + \sum_j e_{s,j} + \sum_j e_{\cap,j} > a_1,\]
then $x_{1,1} m_e/x_{\cap,p} \in (I_d)$, for the all the left-hand sides 
of the inequalities but the first stay the same, and the first inequality for the new monomial is 
\[ \sum_j e_{2,j} + \cdots + \sum_j e_{s,j} + \sum_j e_{\cap,j} \ge a_1.\] 
(Note that this property is independent of whether we use $x_{1,1}$ or 
some other $x_{1,v}$ with $e_{1,v} < f_{1,v}$.) 

If $\sum e_{2,j} + \cdots + \sum e_{s,j} + \sum e_{\cap,j} \not> a_1$, then
\begin{equation} \label{e:int1}
\sum_j e_{2,j} + \cdots + \sum_j e_{s,j} + \sum_j e_{\cap,j} = a_1 \le 
\sum_j f_{2,j} + \cdots + \sum_j f_{s,j} + \sum_j f_{\cap,j}
\end{equation}
since $m_f \in I$. If $e_{2,j} \ge f_{2,j}, \dots,$ 
$e_{s,j} \ge f_{s,j}$ for all $j$, since $\sum e_{\cap,j} > \sum 
f_{\cap,j}$, we have contradicted (\ref{e:int1}). Therefore, without loss of 
generality, we may assume that some $e_{2,j} < f_{2,j}$.

We proceed by induction. Suppose that we have either found an $e_{i,j} < 
f_{i,j}$ such that $x_{i,j} m_e / x_{\cap,p} \in (I_d)$ for some $i \leq 
t-1$,
or for $r=1,\dots,t-1$, we have
\begin{equation}\label{e:intinduct}
\sum_{\stackrel{i=1}{i \not = r}}^s \sum_{j=1}^{b_i} e_{i,j} + 
\sum_{j=1}^{b_{\cap}} e_{\cap,j} = a_r \le \sum_{\stackrel{i=1}{i \not = 
r}}^s \sum_{j=1}^{b_i} f_{i,j} + \sum_{j=1}^{b_{\cap}} f_{\cap,j}.
\end{equation}
(That is, the double sum is the sum of all $e_{i,j}$ with $i \not = r$.)

As part of the induction hypothesis, we may assume that there exists 
$e_{t,j} < f_{t,j}$. If $x_{t,j} m_e / x_{\cap,p} \in (I_d)$, we are done; 
otherwise,
\begin{equation}\label{e:intt}
\sum_{\stackrel{i=1}{i \not = t}}^s \sum_{j=1}^{b_i} e_{i,j} + 
\sum_{j=1}^{b_{\cap}} e_{\cap,j} = a_t \le \sum_{\stackrel{i=1}{i \not = 
t}}^s \sum_{j=1}^{b_i} f_{i,j} + \sum_{j=1}^{b_{\cap}} f_{\cap,j}.
\end{equation}
Summing (\ref{e:intt}) and the inequalities (\ref{e:intinduct}) for all 
$r=1,\dots,t-1$, we obtain 
\begin{eqnarray}\label{intfinalind}
\left ( t-1 \right ) \sum_{i=1}^t \sum_j e_{i,j} + t \sum_{i=t+1}^s \sum_j 
e_{i,j} + t \sum_j e_{\cap,j} \le \\ \left ( t-1 \right ) \sum_{i=1}^t 
\sum_j f_{i,j} + t \sum_{i=t+1}^s \sum_j f_{i,j} + t \sum_j f_{\cap,j}. 
\nonumber
\end{eqnarray}

Subtracting $(t-1)\deg m_e=(t-1)\deg m_f$, we are left with
\[ \sum_{i=t+1}^s \sum_j e_{i,j} + \sum_j e_{\cap,j} \le \sum_{i=t+1}^s 
\sum_j f_{i,j} +  \sum_j f_{\cap,j}.\]
If $e_{t+1,j} \ge f_{t+1,j},\dots,$ $e_{s,j} \ge f_{s,j}$ for all $j$, 
then we have a contradiction since $\sum e_{\cap,j} > \sum f_{\cap,j}$. 
Hence we may assume without loss of generality that some $e_{t+1,j} < 
f_{t+1,j}$.

Therefore either we find some $e_{i,j} < f_{i,j}$, with $1 \le i \le s$, 
such that $x_{i,j} m_e/x_{\cap,p} \in (I_d)$, or else the exchange property is not 
true, and (\ref{e:intinduct}) holds for all $r=1,\dots,s$. In the latter 
case, summing all $s$ inequalities of the form in (\ref{e:intinduct}), we 
have 
\[\left ( s-1 \right ) \sum_{i=1}^s \sum_j e_{i,j} + s \sum_j e_{\cap,j} 
\le \left ( s-1 \right ) \sum_{i=1}^s \sum_j f_{i,j} + s \sum_j 
f_{\cap,j}. \]
If we subtract $(s-1)\deg m_e = (s-1) \deg m_f$ from both sides, we have 
\[ \sum_j e_{\cap,j} \le \sum_j f_{\cap,j}.\]
But this contradicts our assumption that
\[ \sum_j e_{\cap,j} > \sum_j f_{\cap,j}.\]
Hence there exists some $e_{i,j} < f_{i,j}$ such that $x_{i,j} 
m_e/x_{\cap,p} \in (I_d)$, and the exchange condition holds.

The second case to consider is when some $e_{i,j} > f_{i,j}$ for some $1 
\le i \le s$. Without loss of generality, assume that $e_{1,1} > f_{1,1}$. 
If there exists $e_{1,j} < f_{1,j}$, then $x_{1,j} m_e/x_{1,1} \in (I_d)$. 
Otherwise, we have $e_{1,j} \ge f_{1,j}$ for all $j$, and $\sum e_{1,j} > 
\sum f_{1,j}$. Additionally, note that if any $e_{\cap,j} < f_{\cap,j}$, 
then $m'=x_{\cap,j} m_e/x_{1,1} \in (I_d)$ since the degrees of $m_e$ and 
$m'$ in the $J_2,\dots,J_s$ variables are the same, and the degree of $m'$ 
in the $J_1$ variables is one higher than that of $m_e$. Therefore we may 
also assume that $e_{\cap,j} \ge f_{\cap,j}$ for all $j$.

Since $\deg m_e=\deg m_f$, there exists some $e_{i,j} < f_{i,j}$, and we 
may assume that $e_{2,1} < f_{2,1}$. If $x_{2,1} m_e/x_{1,1} \in (I_d)$, 
we are done; otherwise,
\[ \sum_j e_{1,j} + \sum_j e_{3,j} + \cdots + \sum_j e_{s,j} + \sum_j 
e_{\cap,j} = a_2 \le  \sum_j f_{1,j} + \sum_j f_{3,j} + \cdots + \sum_j 
f_{s,j} + \sum_j f_{\cap,j}.\]
If $e_{3,j} \ge f_{3,j}, \dots,$ $e_{s,j} \ge f_{s,j}$ for all $j$, then 
we have a contradiction since $\sum e_{1,j} > \sum f_{1,j}$ and $\sum 
e_{\cap,j} \ge \sum f_{\cap,j}$. Therefore we may assume without loss of 
generality that $e_{3,1} < f_{3,1}$.

Continuing in this way, we apply an almost identical induction argument as 
in the previous case except that now $\sum e_{1,j} > \sum f_{1,j}$ and 
$\sum e_{\cap,j} \ge \sum f_{\cap,j}$. Either we find some $e_{i,j} < 
f_{i,j}$ with $2 \le i \le s-1$ such that $x_{i,j} m_e/x_{1,1} \in (I_d)$, 
or (\ref{e:intinduct}) holds for all $r=2,\dots,s-1$, and there exists some 
$e_{s,j} < f_{s,j}$. If $x_{s,j} m_e/x_{1,1} \in (I_d)$, we are done. 
Otherwise, (\ref{e:intinduct}) also holds for $r=s$, and the exchange 
property fails. Summing the inequalities of (\ref{e:intinduct}) for all 
$r=2,\dots,s$, we obtain
\[\left ( s-1 \right ) \sum_j e_{1,j} + \left ( s-2 \right ) \sum_{i=2}^s 
\sum_j e_{i,j} + \left ( s-1 \right ) \sum_j e_{\cap,j} \le \] \[ \left ( 
s-1 \right ) \sum_j f_{1,j} + \left ( s-2 \right ) \sum_{i=2}^s \sum_j 
f_{i,j} + \left ( s-1 \right ) \sum_j f_{\cap,j}. \]
Now, subtracting $(s-2)\deg m_e = (s-2)\deg m_f$, we have 
\[ \sum_j e_{1,j} + \sum_j e_{\cap,j} \le \sum_j f_{1,j} + \sum_j 
f_{\cap,j}.\]
But we are assuming that $\sum e_{1,j} > \sum f_{1,j}$, and $\sum 
e_{\cap,j} \ge \sum f_{\cap,j}$, so we have a contradiction, and thus the 
exchange property must hold.
\end{proof}

As a consequence of Theorem~\ref{t.maintool}, the intersection
of any two Veronese ideals must be componentwise linear.

\begin{corollary}\label{c.2CWL}
Let $J,K \subseteq [n]$ and let $a$ and $b$ be positive integers. Then
$I = \mf_J^a \cap \mf_K^b \subset R=k[x_1,\dots,x_n]$ is componentwise linear.
\end{corollary}

\begin{proof}
If $J \cup K = [n]$, then we are done by Theorem~\ref{t.maintool}. If not, we may reindex the variables so that 
$J \cup K = [m]$, and $[m] \cup \{m+1,\dots,n\}=[n]$.  Then, by
Theorem~\ref{t.maintool}, $I$ is componentwise linear in $k[x_1,\ldots,x_m]$ and so 
Lemma~\ref{l.extendedring} gives the result. 
\end{proof}

\begin{remark} \label{r.notpoly}
It is not true that all ideals $I=\mf_{J_1}^{a_1} \cap \mf_{J_2}^{a_2}$ have $(I_d)$ polymatroidal for all $d$.
 For example, let 
\[ I=\mf_{J_1}^{a_1} \cap \mf_{J_2}^{a_2} = (x_1,x_2,x_3)^2 \cap (x_2,x_3,x_5)^2 \subset R=k[x_1,\dots,x_5] .\]
 Note that $x_4$ does not appear in either of 
the two components. Both $m_e=x_3^2x_4$ and $m_f=x_1x_3x_5$ are in $(I_3)$. The power on $x_3$ is greater in
 $m_e$ than it is in $m_f$, and the powers on $x_1$ and $x_5$ are larger in $m_f$ than in $m_e$. But 
$x_1(x_3^2x_4/x_3)=x_1x_3x_4 \not \in (I_3)$, and $x_5(x_3^2x_4/x_3)=x_3x_4x_5 \not \in (I_3)$. Therefore 
$(I_3)$ is not polymatroidal. The proof of Theorem~\ref{t.maintool} breaks down because $4$ is missing 
from both $J_1$ and $J_2$, and hence we cannot partition the variables the way we did in that argument;
 the $x_4$ exponents would be double-counted, causing problems when we subtract a multiple of $\deg m_e=\deg m_f$.
\end{remark}

\section{The intersection of three Veronese ideals} \label{s.threeideals}

We will show that intersection of any three Veronese ideals is always
componentwise linear.  Throughout this section,  we write 
$\G(I)$ to denote the set of minimal generators of a monomial ideal $I$.

We begin with an observation.  Suppose $J,K \subseteq [n]$, but 
$J \cup K \subsetneq [n]$. Let $H = [n] \backslash
(J \cup K)$, and after relabeling, we can assume $H = \{r+1,\ldots,n\}$
and $J \cup K = [r]$.
Let $a,b$ be any positive integers, and let $\alpha$ be the smallest
degree of a nonzero element in $(\mf_J^a \cap \mf_K^b)$.  If 
\[B' = (\mf_J^a \cap \mf_K^b) \cap k[x_i ~|~ i \in J \cup K] = (\mf_J^a 
\cap \mf_K^b) \cap 
k[x_1,\ldots,x_r], \] then
the ideal $B = ((\mf_J^a \cap \mf_K^b)_{\alpha+d})$
as an ideal of $R =k[x_1,\ldots,x_n]$ has the decomposition
\[B = (B'_{\alpha+d}) + (B'_{\alpha+d-1})\mf_H + (B'_{\alpha+d-2})\mf_H^2 
+ \cdots + (B'_{\alpha})\mf_H^d, \]
where $(B'_i)$ denotes the ideal generated by elements of degree $i$ of 
$B'$ but viewed as an ideal of $R$.

Order the elements of $\G(B)$ as follows: Order the generators 
of $(B'_{\alpha+d})$ with respect to the  ascending reverse-lex 
ordering.  Then add the generators of $(B'_{\alpha+d-1})\mf_H$ in ascending 
reverse-lex order.  We thus add all the monomials
divisible by $x_n$ first.  Continue by adding the 
generators $(B'_{\alpha+d+2})\mf_H^2$
in ascending reverse-lex order, and so on.

\begin{lemma}\label{two ideal lemma}
Using the above notation, the ideal
\[B = (B'_{\alpha+d}) + (B'_{\alpha+d-1})\mf_H + (B'_{\alpha+d-2})\mf_H^2 
+ \cdots + (B'_{\alpha})\mf_H^d\]
has linear quotients with respect to the order prescribed above.
\end{lemma}

\begin{proof}
Let $M_i$ be the $i$-th element of $\G(B)$, $i \geq 2$,
with respect to our ordering.  First,
suppose that $M_i \in (B'_{\alpha+d})$.  We wish to calculate $(M_1,\ldots,M_{i-1}):M_i$,
where $(M_1,\ldots,M_{i-1})$ is the ideal generated by all monomials in 
$\G(B)$ smaller
than $M_i$ with respect to our ordering.  Note that in this case, each 
$M_j$ is in $(B'_{\alpha+d})$,
so each $M_j$ is in $S = k[x_1,\ldots,x_r]$.  As an ideal of $S$,
the ideal $(B'_{\alpha+d})$ is polymatroidal by Theorem~\ref{t.maintool}.  So, $(B'_{\alpha+d})$ 
has linear quotients
with respect to the ascending reverse-lex order
by  Proposition~\ref{p.poly-asc}.  Thus, as an ideal of $S$,
$(M_1,\ldots,M_{i-1}):M_i = (x_{i_1},\ldots,x_{i_j})$ for some subset
$\{i_1,\ldots,i_j\} \subseteq J \cup K$.  Because $S \rightarrow R$ is a flat ring homomorphism, by \cite[Theorem 7.4(iii)]{Mat}, $(M_1,\ldots,M_{i-1})R:M_iR = (x_{i_1},\ldots,x_{i_j})R$. (We note that if $I$ is an ideal of S, then we will sometimes abuse notation and write $I$ to denote an ideal of $R$, where we really mean $IR$ using the flat homomorphism $S \rightarrow R$.)

So, suppose now that $M_i \in (B'_{\alpha+d-s})\mf_H^s$ for some $s = 
1,\ldots,d$.  Let
\[I = (B'_{\alpha+d}) + \cdots +
(B'_{\alpha+d-s+1})\mf_H^{s-1} + (M \in \mathcal{G}((B'_{\alpha+d-s})\mf_H^s) ~|~ M < 
M_i),\]
where $M < M_i$ with respect to the reverse-lex ordering.  

Since $M_i \in (B'_{\alpha+d-s})\mf_H^s$,
$M_i = M_{i1}M_{i2}$ with $M_{i1} \in (B'_{\alpha+d-s})$ and $M_{i2} \in 
\mf_H^s$.  If we multiply $M_i$
by any $x_l$ with $l \in J \cup K$, then $x_lM_{i1} \in 
(B'_{\alpha+d-s+1})$.  Since $M_{i2} \in
\mf_H^s \subseteq \mf_H^{s-1}$, we have $x_lM_i \in 
(B'_{\alpha+d-s+1})\mf_H^{s-1} \subseteq I$.  Thus
$\mf_{J\cup K} \subseteq I:M_i$. 

Because $M_{i2} \in \mf_H^s$, we have
\[M_{i2} = x_{r+1}^{c_{r+1}}\cdots x_n^{c_n} ~~\mbox{with $c_{r+1} + 
\cdots + c_n = s$}.\]
Let $l$ be the smallest integer in $\{r+1,\ldots,n\}$ such that $c_l > 
0$.  Then
$x_eM_i \in I$ for $e = l+1,\ldots,n$.  To see this note that
\[x_eM_i = M_{i1}x_eM_{i2} = M_{i1}x_{l}^{c_{l}} \cdots x_e^{c_e+1} \cdots 
x_n^{c_n} = 
x_l M_{i1}x_{l}^{c_l-1}\cdots x_e^{c_e+1} \cdots x_n^{c_n}.\]
But then $M_{i1}x_{l}^{c_l-1}\cdots x_e^{c_e+1} \cdots x_n^{c_n} < M_i$ 
since $c_e+1 > c_e$.
So 
\[x_eM_i = x_l M_{i1}x_{l}^{c_l-1}\cdots x_e^{c_e+1} \cdots x_n^{c_n} \in (M \in \G((B'_{\alpha+d-s})\mf_H^s) ~|~ M < M_i) \subseteq I,\]
from which we deduce that $(x_{l+1},\ldots,x_n) \subseteq I:M_i$.

Let $m$ now be any monomial not in $\mf_{J\cup K} + (x_{l+1},\ldots,x_n)$ 
and suppose
$mM_i \in I$.  The monomial $m$ can only be divisible
by the variables $x_{r+1},\ldots,x_l$; suppose $\deg m = z$.
Since $mM_i \in I$, there exists a monomial $M_j \in I$ such that $mM_i = 
m'M_j$
for some monomial $m'$.  Since $M_i \in \mf_H^{s}$, $mM_i \in 
\mf_H^{s+z}$.
If $M_j \in (B'_{\alpha+d-i})\mf_H^i$ for some $i < s$, then $m'M_j$
cannot be in $\mf_H^{s+z}$ since in $m'M_j$, the exponents of 
$x_{r+1},\ldots,x_n$ can
add up to at most $i+z$.  Thus, we must have $M_j \in 
(B_{\alpha+d-s})\mf_H^s$, and so
$M_j < M_i$ with respect to the reverse-lex ordering.

If we write $M_j = M_{j1}M_{j2}$ with
$M_{j1} \in (B_{\alpha+d-s})$ and $M_{j2} \in \mf_H^s$, then since $m$ is a 
monomial
in the variables $x_{r+1},\ldots,x_l$ only, we must 
have $M_{j1} = M_{i1}$
and $M_{j2} < M_{i2}$.  If $M_{j2} = x_{r+1}^{f_{r+1}}\cdots x_n^{f_n}$,
there must be some $e$ such that $f_e > c_e$ but $f_{e+1} = c_{e+1}, \ldots 
f_{n} = c_n$.
Furthermore, since 
$M_{i2} = x_l^{c_l} \cdots x_n^{c_n}$, we must have $l < e \leq n$.  
Indeed, if $e \le l$, then
\[ s = f_{r+1} + \cdots + f_n \ge f_e + \cdots f_n > c_e + \cdots + c_n = c_l + \cdots + c_n  = s.\]
Thus, for  $mM_i = m'M_j$ to be true, both sides must be divisible by 
$x_e^{f_e}$.  But since
$M_i$ is not divisible by $x_e^v$ with $v > c_e$, we must have $m$ 
divisible by $x_e$.
But this contradicts the fact that $m$ is not in the ideal $\mf_{J\cup K} 
+ (x_{l+1},\ldots,x_n)$.
We then arrive at the conclusion
\[ I:M_i = \mf_{J\cup K} + (x_{l+1},\ldots,x_n).\]
So, $B$ has linear quotients.
\end{proof}

\begin{remark} \label{r.anotherpf}
Lemma~\ref{two ideal lemma} gives a second proof that $\mf_J^a \cap \mf_K^b$ is 
componentwise linear for all $J$, $K$, $a$, and $b$.
\end{remark}

We thank the referee for suggestions that significantly simplified the proof of the following theorem.

\begin{theorem} \label{t.threeideals}
Let $J_1,J_2,J_3 \subseteq [n]$ be three sets, and let $a_1,a_2,a_3$ be
three positive integers.  Then $\mf_{J_1}^{a_1} \cap \mf_{J_2}^{a_2} \cap 
\mf_{J_3}^{a_3}$ is componentwise linear.
\end{theorem} 

\begin{proof} 
If $J_i = J_j$ for some $i$ and $j$, then $\mf_{J_i}^{a_i} \cap \mf_{J_j}^{a_j} = 
\mf_{J_i}^{\max\{a_i,a_j\}}$, and thus we are in the case of Corollary \ref{c.2CWL}.  So,
we may assume that all the $J_i$'s are distinct.  Next
we may assume by Lemma~\ref{l.extendedring} that 
$J_1 \cup J_2 \cup J_3 = [n]$. If the hypotheses of Theorem~\ref{t.maintool} are 
satisfied, we done.  So we may further assume that there exists a pair of sets 
$J_i$ and $J_j$ such that $J_i \cup J_j \subsetneq [n]$.  

For ease of exposition, we shall use $J,K,L$
for $J_1,J_2,J_3$, we shall use $a,b,c$ for $a_1,a_2,a_3$,
and we shall assume that $J \cup K \subsetneq [n]$ and that $J$, $K$, and $L$ are distinct. 
After relabeling,
we can also assume that $L =  \{t,t+1,\ldots,n\}$.
We also set  $H_1 = L \cap (J \cup K)$ and $H_2 = L \backslash (J\cup K)$.  
After relabeling
again, we may further assume that $H_1 = \{t,\ldots,r\}$ and $H_2 = 
\{r+1,\ldots,n\}$.

Let $\alpha$ be the smallest degree of a nonzero element in $\mf_J^a \cap 
\mf_K^b$.  Because
$(\mf_J^a \cap \mf_K^b \cap \mf_L^c)_e \subseteq (\mf_J^a \cap \mf_K^b)_e$ 
for all $e$, $(\mf_J^a \cap \mf^b_K \cap \mf_L^c)_e = (0)$ if $e < \alpha$, and thus 
has a linear resolution.

Now fix a $d \geq 0$, and set  
$A = ((\mf_J^a \cap \mf_K^b \cap \mf_L^c)_{\alpha+d})$ and $B = ((\mf_J^a 
\cap \mf_K^b)_{\alpha+d})$.
We shall show that $A$ has linear quotients, and hence $A$ 
has a linear resolution.
It will then follow that   $\mf_{J}^{a} \cap \mf_{K}^{b} \cap \mf_{L}^{c}$
is componentwise linear.

Set \[
B' = (\mf_J^a \cap \mf_K^b) \cap k[x_i ~|~ i \in J \cup K] = 
(\mf_J^a \cap \mf_K^b) \cap k[x_1,\ldots,x_r].\] 
Note that the ideal $B'$ has the same generators as $\mf_J^a \cap \mf_K^b$, but we are now 
considering $B'$ as
an ideal in a smaller ring.
 The ideal $B$ then has
the following decomposition:
\[B = (B'_{\alpha+d}) + (B'_{\alpha+d-1})\mf_{H_2} + 
(B'_{\alpha+d-2})\mf_{H_2}^2 + \cdots + (B'_{\alpha})\mf_{H_2}^{d}\]
where by $(B'_i)$ we mean the ideal generated by the degree $i$ part of 
$B'$ in $k[x_1,\ldots,x_r]$, but
considering the ideal as an ideal of $R$. 

Since $A \subseteq B$, each generator of $A$ belongs to some $(B'_{\alpha+d 
-i})\mf_{H_2}^i$
for some $i = 0,\ldots,d$.  Set 
\[A_i = \{M ~|~ M \in \G(A) ~\mbox{and}~ M \in 
(B'_{\alpha+d-i})\mf_{H_2}^i\} 
~\mbox{for each $i =0,\ldots,d$}.\]
So $\G(A) = A_0 \cup A_1 \cup \cdots \cup A_d$.

Order the elements of $\G(A)$ as follows:  Begin by adding the 
elements of $A_0$
in ascending reverse-lex order.  Then, add the elements of $A_1$, after all 
the elements in $A_0$, in ascending reverse-lex order.  We then add 
the elements of $A_2$
in ascending reverse-lex order, and so on.  The ordering could also be 
described as follows: Write out
the generators of $B$ in the same order as  in Lemma \ref{two ideal 
lemma}.  Then simply remove
any element of $\G(B)$ that is not in $\G(A)$.
We will show that $A$ has linear quotients
with respect to this ordering.

Let $M_i$ be the $i$-th monomial of $\G(A)$ with respect to our ordering 
with $i \geq 2$.
Set $I = (M_1,\ldots,M_{i-1})$, the ideal generated by all the monomials 
in $\G(A)$
smaller than $M_i$ with respect our ordering.  Furthermore, let $D = (M 
\in \G(B) ~|~ M < M_i)$ with respect
to our ordering.  Since $I \subseteq D$, by Lemma \ref{two ideal lemma} we
have 
\[I:M_i \subseteq D:M_i = (x_{i_1},\ldots,x_{i_j})\]
because $B$ has linear quotients with respect to this order.

If $M_i \in A_0$, then we will show that $I:M_i = D:M_i$.  Take any 
$x_e \in \{x_{i_1},\ldots,
x_{i_j}\}$.  Then $x_eM_i = x_fM_j$ for some $M_j \in D$.  We wish to show
that $M_j$ is in $I$.  Suppose $M_j \not\in I$, i.e., $M_j \not\in A$.  
Then 
\[M_j = x_1^{d_1} \cdots x_t^{d_t} \cdots x_n^{d_n} ~~\mbox{with $d_t + 
\cdots + d_n < c$.}\]
On the other hand, $M_i \in A$, so
\[M_i = x_1^{c_1} \cdots x_t^{c_t} \cdots x_n^{c_n} ~~\mbox{with $c_t + 
\cdots + c_n \geq c$.}\]
Since $M_i \in A_0 \subseteq (B'_{\alpha+d})$, we have $M_j \in 
(B'_{\alpha+d})$ as well, and $M_j < M_i$ with respect to reverse-lex order.  Hence there must exist 
some $l$
such that $d_l > c_l$, but $d_{l+1} = c_{l+1}, \ldots, d_n = c_n$.  Moreover, $l 
\in \{t,\ldots,n\}$,
for if $l < t$,  then 
\[c > d_t + \cdots + d_n = c_t + \cdots + c_n \geq c.\]
Thus, for $x_eM_i = x_fM_j$ to be true, $x_e = x_l$, since the exponent
of $x_l$ is higher in $M_j$ than in $M_i$.  So the exponents of 
$x_t,\ldots,x_n$ in $x_eM_i$
must add up at least $c+1$.  However, the exponents of $x_t,\ldots,x_n$ 
in $x_fM_j$
can add up to at most $c$, a contradiction.  So $M_j$ must be in $A$ and thus is in $I$.

Suppose now that $M_i \in A_s$ with $s \in \{1,\ldots,d\}$.  So $M_i \in 
(B'_{\alpha+d-s})\mf_{H_2}^s$.
Write $M_i$ as $M_i = M_{i1}M_{i2}$ with $M_{i1} \in (B'_{\alpha+d-s})$ 
and 
$M_{i2} \in \mf_{H_2}^s.$  Then 
\[M_{i2} = x_{r+1}^{c_{r+1}} \cdots x_n^{c_n} ~~\mbox{with $c_{r+1} + 
\cdots + c_n = s$.}\]
Let $l$ be the smallest integer in the set $\{r+1,\ldots,n\}$ such that 
$c_l > 0$.  As
shown in the proof of Lemma \ref{two ideal lemma},
\[D:M_i = \mf_{J\cup K} + (x_{l+1},\ldots,x_n).\]
Since $I:M_i \subseteq D:M_i$, the above fact implies that no monomial
of the form $x_{r+1}^{c_{r+1}}\cdots x_l^{c_l}$ can belong to $I:M_i$.

Next, we show that $(x_{l+1},\ldots,x_n) \subseteq I:M_i$.  Let $x_e \in 
\{x_{l+1},\ldots,x_n\}$.
Then
\[x_eM_i = M_{i1}x_l^{c_l}\cdots x_e^{c_e+1} \cdots x_n^{c_n} = 
x_lM_{i1}x_l^{c_l-1}\cdots
x_e^{c_e+1} \cdots x_n^{c_n}, \]
and $M_{i1}x_l^{c_l-1}\cdots x_e^{c_e+1} \cdots x_n^{c_n} \in (B'_{\alpha+d-s})\mf_{H_2}^s$.  Also, it is clear that 
\[M_{i1}x_l^{c_l-1}\cdots x_e^{c_e+1} \cdots x_n^{c_n} \in A_s\] since the 
exponents of $x_t,\ldots,x_n$ still add up
to at least $c$.  Now 
\[M_j = M_{i1}x_l^{c_l-1}\cdots
x_e^{c_e+1} \cdots x_n^{c_n} < M_i\] with respect to the reverse-lex 
order, so $M_j \in I$,
and hence $x_e \in I:M_i$.

Now suppose that $x_e \in \{x_1,\ldots,x_r\} = \{x_i ~|~ i \in J \cup 
K\}$.  Since 
$x_e \in D:M_i$, we have that $x_eM_i$ is divisible by some monomial $M 
\in D$ with
$M$ less than $M_i$ with respect to our ordering.
The monomial $M$ may or may not be in $I$.  We thus
partition $J \cup K$ into the following two sets:
\begin{eqnarray*}
P_1 &= &\{e \in J \cup K ~|~ x_eM_i ~\mbox{is divisible by some $M \in 
D$ with $M \in I$}\} \\
P_2 & =&\{e \in J \cup K ~|~ \mbox{for every $M \in D$ such that 
$M|x_eM_i$, $M \not\in I$}\}.
\end{eqnarray*}
It follows immediately that if $e \in P_1$, then $x_e \in I:M_i$, so 
$\mf_{P_1} 
\subseteq I:M_i$.

We will now show (through many steps) that if $m$ is any monomial in the variables
$\{x_e ~|~ e \in P_2\}$, then $m \not\in I:M_i$.  It then follows that 
\[I:M_i = \mf_{P_1} + (x_{l+1},\ldots,x_n)\]
so that $I$ has linear quotients.

Suppose that $e \in P_2$.  Then $x_eM_i = x_fM_j$ for some $M_j \in D$,
and also $M_j \not\in I$.  Furthermore, let
\[ M_i = x_1^{c_1}\cdots x_n^{c_n} ~~\text{and}~~ M_j = x_1^{d_1}\cdots 
x_n^{d_n}.\]

We begin with some facts that must be true about $x_e,x_f,M_i$ and $M_j$
in this case.  First, since $M_j \not\in I$, $d_t+\cdots +d_n <c$,
but $M_i \in I$ means that $c_t + \cdots + c_n \geq c$.  Since
$x_fM_j = x_eM_i$, this implies that $x_f \in \{x_t,\ldots,x_n\}$
and $x_e \in \{x_1,\ldots,x_{t-1}\}$.  We also observe that
this must imply that
\[d_t+\cdots+d_n = c-1 ~~\text{and}~~ c_t+\cdots + c_n = c.\]

Second, the monomial $M_j \not\in (B'_{\alpha+d-s})\mf_{H_2}^s$.
Observe that if $M_j \in (B'_{\alpha+d-s})\mf_{H_2}^s$,
then $M_j < M_i$ with respect to the reverse-lex order.  So there
exists some index $p$ such that $d_p > c_p$ but $d_{p+1} = c_{p+1},
\ldots,d_n = c_n$.  Now because $x_eM_i = x_fM_j$, we must
have $x_e = x_p$.  But since $x_p \in \{x_1,\ldots,x_{t-1}\}$ we would have
\[c > d_t+\cdots+d_n = c_t+\cdots+c_n \geq c,\]
which is a contradiction.  Note that this argument
applies to any monomial $M \in D$
with the property that $M|x_eM_i$ but $M \not\in I$.

Now, let $m$ be any monomial in the variables $\{x_e ~|~ e \in P_2\}$, and
suppose that
$mM_i \in I$; that is,
$m \in I:M_i$.  Then there exists a monomial $m' \in R$ and $M \in \G(I)$
such that
$mM_i = m'M$.  If $M = x_1^{b_1} \cdots x_n^{b_n}$, then $b_t + \cdots +
b_n \geq c$
since $M \in I$.  Since $m$ is not divisible by any element of
$\{x_t,\ldots,x_n\}$  (this
follows since any variable in $\{x_e ~|~ e \in P_2$ must be in
$\{x_1,\ldots,x_{t-1}\}$), the
exponents of $x_t,\ldots,x_n$ in $mM_i$ are the same as those of $M_i$,
and
thus, the exponents of $x_t,\ldots,x_n$ in $mM_i$ add up to $c$
since $c_t+\cdots+c_n = c$.  Thus,
any variable that divides $m'$ must also be in $\{x_1,\ldots,x_{t-1}\}$;
otherwise, the exponents of $x_t,\ldots,x_n$ in $m'M$ add up to a number greater
than $c$.

Since $m$ and $m'$ are only divisible by the variables $x_1,\ldots,x_{t-1}$,
we
must therefore
have $b_t = c_t, \ldots, b_n = c_n$.  In particular, $b_{r+1} = c_{r+1},
\ldots, b_n = c_n$.
Thus, if we let $M^{\star} = x_{r+1}^{b_{r+1}} \cdots x_n^{b_n}$, then
$M_i = M'_iM^{\star}$ and $M = M'M^{\star}$.  Because $M^{\star} \in
\mf_{H_2}^s$,
we have $M,M_i \in (B'_{\alpha+d-s})\mf_{H_2}^s$.  It follows that $M'_i,
M' \in (B'_{\alpha+d-s})$.
Since $M$ is in $I$, we have $M < M_i$ with respect to the reverse-lex
ordering.
This implies that $M' < M'_i$ with respect to the reverse-lex
ordering.

Let $D'$ be the ideal of $S = k[x_1,\ldots,x_r]$
generated by all generators of $(B'_{\alpha+d-s})$
less than $M'_i$
with respect to the reverse-lex order.  Since
$(B'_{\alpha+d-s})$ is
polymatroidal in this ring by Theorem \ref{t.maintool},
it has linear quotients with respect to the ascending reverse-lex ordering
(by
 Proposition~\ref{p.poly-asc}).
So
\[ D':M'_i = (x_{i_1},\ldots,x_{i_g}).\]
Now $mM_i = m'M$ implies that $mM'_i = mM'$.  Since $m,m'$ can be viewed
as elements of $S$,
we have $m \in D':M'_i$ since $M' \in D$.  So there exists some $x_e \in
\{x_{i_1},\ldots,x_{i_g}\}$
such that $x_e |m$.  Note that $e \in P_2$.  We thus must have some $M''
\in D'
\subseteq (B'_{\alpha+d-s})$ with $M'' < M'_i$
such that $x_eM'_i = x_fM''$ for some $x_f$.  But then
\[x_eM_i = x_eM'_iM^{\star} = x_fM''M^{\star}.\]
Now $M''M^{\star} \in D$ since $M''M^{\star} < M'_iM^{\star} = M_i$.
We must have that $M''M^{\star} \not\in I$, because
if $M''M^{\star} \in I$, this would imply that $x_e \in P_1$
since $M''M^{\star}|x_eM_i$.
But we also have that $M''M^{\star} \in
(B'_{\alpha+d-s})\mf_{H_2}^s$,
contradicting the fact that  every $M \in D$
with the property  $M|x_eM_i$ but $M \not\in I$ cannot be
in $(B'_{\alpha+d-s})\mf_{H_2}^{s}$.

Thus $m \not\in I:M_i$, and the conclusion follows.
\end{proof}

\begin{remark} \label{r.simplest}
Combining Corollary~\ref{c.2CWL} and Theorem~\ref{t.threeideals}, we conclude that
\[ I=(x_1,x_2) \cap (x_2,x_3) \cap (x_3,x_4) \cap (x_1,x_4) \]
is the simplest intersection of Veronese ideals that is not componentwise linear. It is the ideal of a 
tetrahedral curve; see \cite{MN} and \cite{FMN} for studies of these ideals and their resolutions, 
including a characterization of which curves are componentwise linear. Note that to form an intersection 
of Veronese ideals that is not componentwise linear, by our earlier results, we must have $s \ge 4$. By 
analyzing the possible cases for three variables, it is not hard to see that we must also work in a ring
 with at least four variables: The presence of any ideal $(x_1,x_2,x_3)^a \subset k[x_1,x_2,x_3]$ in the
 intersection is irrelevant to componentwise linearity, and hence one needs only show that
\[ (x_1)^{a_1} \cap (x_2)^{a_2} \cap (x_3)^{a_3} \cap (x_1,x_2)^{a_4} \cap (x_1,x_3)^{a_5} \cap (x_2,x_3)^{a_6}, \]
where the $a_i \ge 0$, can be expressed as
\[ (x_1^{a_1}x_2^{a_2}x_3^{a_3}) \left ( (x_1,x_2)^{\max\{a_4-(a_1+a_2),0\}} \cap 
(x_1,x_3)^{\max\{a_5-(a_1+a_3),0\}} \cap (x_2,x_3)^{\max\{a_6-(a_2+a_3),0\}} \right ). \]
Theorem~\ref{t.threeideals} tells us that this ideal is componentwise linear.
\end{remark}

The proof of Theorem~\ref{t.threeideals} gives some insight into why there are ideals with $s=4$ that
 fail to be componentwise linear. The ideals 
\[B'=(\mf_J^a \cap \mf_K^b) \cap k[x_i ~|~i \in J \cup K]=(\mf_J^a \cap \mf_K^b) \cap k[x_1,\dots,x_r]\] 
play a prominent role in the proof. We use the fact that the ideals $(B'_{\alpha+d-s})$ are polymatroidal
 in $k[x_1,\dots,x_r]$ by Theorem~\ref{t.maintool}, which shows that they have linear quotients with
 respect to ascending reverse-lex order. If, in trying to prove the $s=4$ case, we defined the $B'$ as 
the intersection of three ideals $\mf_J^a$, $\mf_K^b$, and $\mf_L^c$, intersected with the appropriate 
ring, this step would fail without extra hypotheses on $J$, $K$, and $L$.


\section{Resolutions of $\mf_J^a \cap \mf_K^b$} \label{s.twoideals}

In this section we provide a thorough analysis of the graded Betti numbers of ideals of the form  
$I = \mf_{J}^{a} \cap \mf_{K}^{b}$ with $a \geq b \geq 1$. We derive formulas for the Betti numbers of 
these intersections of Veronese ideals that enable us in the next section to recapture 
the formulas of Valla \cite{Va}, Fatabbi and Lorenzini \cite{FL}, and the first author \cite{Fra} for the 
graded Betti numbers of two fat points in $\pr^n$.  In fact, we can extend
their results to compute the $\N$-graded Betti numbers of two fat points in the multiprojective space 
$\pr^{n_1} \times \cdots \times \pr^{n_r}$. 

To compute the graded Betti numbers of $I = \mf_J^a \cap \mf_K^b$, we generalize the approach
given by the first author in \cite{Fra}.  Our proof hinges on the fact that
$I$ is an example of a splittable monomial ideal.
As in the previous section, for a monomial ideal $I$ we let $\G(I)$ denote the unique set of minimal generators of $I$.
 
\begin{definition}[see \cite{EK}]\label{defn: split}
A monomial ideal $I$ is {\bf splittable} if $I$ is the sum
of two nonzero monomial ideals $J$ and $K$, that is, $I = J+K$, such
that
\begin{enumerate}
\item $\G(I)$ is the disjoint union of $\G(J)$ and $\G(K)$.
\item there is a {\bf splitting function}
\begin{eqnarray*}
\G(J\cap K) &\rightarrow &\G(J) \times \G(K) \\
w & \mapsto & (\phi(w),\psi(w))
\end{eqnarray*}
satisfying
\begin{enumerate}
\item for all $w \in \G(J \cap K), ~~ w = \lcm(\phi(w),\psi(w))$.
\item for every subset $S \subset \G(J \cap K)$, both
$\lcm(\phi(S))$ and $\lcm(\psi(S))$
strictly divide $\lcm(S)$.
\end{enumerate}
\end{enumerate}
If $J$ and $K$ satisfy the above properties, then
we shall say $I = J + K$ is a {\bf splitting} of $I$.
\end{definition}
 
When $I = J + K$ is a splitting of a monomial
ideal $I$, then there is a relation between  $\beta_{i,j}(I)$
and the graded Betti numbers of the smaller ideals. 
 
\begin{theorem}[Eliahou-Kervaire \cite{EK}, Fatabbi \cite{Fa}]
\label{EKF}
Suppose $I$ is a
splittable monomial ideal with splitting $I = J+K$.  Then
for all $i, j \geq 0$,
\[\beta_{i,j}(I) = \beta_{i,j}(J) + \beta_{i,j}(K) +
\beta_{i-1,j}(J\cap K).\]
\end{theorem}

The following lemma (for a proof see \cite[Lemma 1.5]{HT} or \cite[Lemma 2.3]{Fra}) will allows us to determine when a resolution built via a mapping cone construction is in fact minimal.

\begin{lemma}\label{mapping cone lemma}
Let $I \subseteq R = k[x_1,\ldots,x_n]$ be a homogeneous ideal with the regularity
of $R/I$ at most $d-1$.  Let $m$ be a monomial of degree $d$ not in $I$ such that
$I:m = \mf_J$ for some $J \subseteq [n]$.  Then the mapping 
cone resolution of $R/(I,m)$ is minimal.
\end{lemma}

With these tools we can now turn to the graded Betti numbers of 
$I = \mf_{J}^{a} \cap \mf_{K}^{b}$.  The resolution depends upon
how the two subsets $J,K \subseteq [n]$ intersect. There are four possible cases, as listed below,
and we shall deal with each case separately.
\vspace{.25cm}

\noindent
{\bf Case 1: $J \cap K = \emptyset$.}  

If $J \cap K = \emptyset$, then 
$I = \mf_J^{a} \cap \mf_K^b = \mf_J^a\mf_K^b$.  The resolution of $I$ is then
a corollary of Theorem \ref{betti numbers many J}.  For completeness we record the formula here:
\[\beta_{i,i+a+b}(I) = \sum_{i_1+i_2 = i} \binom{|J| + a -1}{a+i_1}\binom{a+i_1 -1}{i_1}
\binom{|K| + b -1}{b+i_2}\binom{b+i_2 -1}{i_2}\]
and $\beta_{i,j}(I) = 0$ for all other $i,j \geq 0$.
\vspace{.25cm}

\noindent
{\bf Case 2: $J \backslash K = \emptyset$} (i.e., $J \subseteq K$).

In this case $I = \mf_J^a \cap \mf_K^b = \mf_J^a.$
By Lemma \ref{betti numbers one J}, the resolution of $I$ is then given by
\[\beta_{i,i+a}(I) = \binom{|J|+a-1}{a+i}\binom{a+i-1}{i}
~\mbox{and $\beta_{i,j}(I) = 0$ otherwise.}\]
\vspace{.25cm}

\noindent
{\bf Case 3: $K \backslash J = \emptyset$} (i.e., $K \subseteq J$).  

Set $A = J \backslash K$,
and let $\mf_A$ denote the corresponding ideal.  In this situation 
$I = \mf_J^a \cap \mf_K^b = \mf^a_K + \mf_A\mf_K^{a-1} + \cdots + \mf_A^{a-b}\mf_K^b.$
We will postpone describing  $\beta_{i,j}(I)$ in this case since these numbers  will
be a byproduct of our work in the final case.
\vspace{.25cm}

\noindent
{\bf Case 4: $J \cap K, J \backslash K, K \backslash J \neq \emptyset$.}

Set $A = J \backslash K, ~~ B = K \backslash J$  and $ C  = J \cap K.$
Let $\mf_A, \mf_B$ and $\mf_C$ be the corresponding monomial ideals. 

\begin{notation}
Since the generators of $\mf_A, \mf_B$ and $\mf_C$ are disjoint subsets of indeterminates of $R$,
for ease of exposition we write
$\mf_A =\langle x_1,\ldots,x_{t_1} \rangle,
\mf_B = \langle y_1,\ldots,y_{t_2} \rangle,$ and $\mf_C = \langle z_1,\ldots,z_{t_3}\rangle.$
\end{notation}

With this notation, we set
\begin{eqnarray*}
U &= &\mf_C^{a} + \mf_A\mf_C^{a-1}+ \cdots + \mf_A^{a-b}\mf_C^{b} \\
V & =&\mf_B\mf_A^{a-b+1}\mf_C^{b-1} + \mf_B^2\mf_A^{a-b+2}\mf_C^{b-2} + \cdots + \mf_B^{b}\mf_A^{a}.
\end{eqnarray*}
To find the graded Betti numbers of $I$, we will exploit the fact that $U$ and $V$ form
a splitting of $I$ (as we prove below).

\begin{theorem}
Suppose that  $J,K$ are subsets of $[n]$ that belong to Case 4, and let $U$ and $V$
be as above.  Then $I =\mf_J^{a} \cap \mf_K^{b}$ is a splittable ideal with splitting
$I = U + V$.
\end{theorem}

\begin{proof}
It is easy to check that $I = U+V$. The containment $U+V \subseteq I$ follows directly from the definitions of $U$ and $V$, and the other containment is a consequence of the fact that $\mf_A$, $\mf_B$, and $\mf_C$ are generated by disjoint monomials.

The definition of $U$ and $V$ gives $\G(U) \cap \G(V) = \emptyset$.  To show
that $U$ and $V$ is a splitting, we first make the observation that
\[U \cap V = \mf_B\mf_A^{a-b+1}\mf_C^b,\]
and hence
\[\G(U \cap V) = \{y_im_1m_2 ~|~ y_i \in \G(\mf_B), ~~ m_1 \in \G(\mf_A^{a-b+1}), ~~ m_2 \in \G(\mf_C^b)
\}.\]
We define our splitting function as follows:
\begin{eqnarray*}
\G(U \cap V) & \rightarrow & \G(U) \times \G(V) \\
m = y_im_1m_2 & \mapsto & (\phi(m),\psi(m)) = ((m_1/x_{\max(m_1)})m_2,y_im_1(m_2/z_{\max(m_2)}))
\end{eqnarray*}
where $\max(m_1) = \max\{i ~|~ x_i|m_1 \}$ and $\max(m_2) = \max\{i ~|~ z_i|m_2\}$.  
It is immediate that $\lcm(\phi(m),\psi(m)) = m$, so our splitting function satisfies the first
condition. 

To verify the second condition, let $S \subseteq \G(U \cap V)$.  It is
straightforward to check that both $\lcm(\phi(S))$ and $\lcm(\psi(S))$ divide $\lcm(S)$.
Moreover, $\lcm(\phi(S))$ strictly divides $\lcm(S)$ since $\lcm(S)$ is divisible
by some $y_{\ell}$, but $\lcm(\phi(S))$ is not.  To see that $\lcm(\psi(S))$ strictly
divides $\lcm(S)$, let $m =y_im_1m_2\in S$ be the monomial with largest ${\max(m_2)}$,
and among all monomials $m' \in S$ divisible by $z_{\max(m_2)}$, the power of $z_{\max(m_2)}$
in $m$, say $d$, is the largest. Hence $z_{\max(m_2)}^d | \lcm(S)$, but $z_{\max(m_2)}^d$
does not divide $\lcm(\psi(S))$.  This implies that $\lcm(\psi(S))$ strictly divides $\lcm(S)$.

So, $I = U+V$ is a splitting of $I$.
\end{proof}

Since $U$ and $V$ is a splitting of $I$, by Theorem \ref{EKF} we only need to compute
the graded Betti numbers of $U, V$, and $U \cap V$.  As noted within the previous
proof
\[U \cap V = \mf_B\mf_A^{a-b+1}\mf_C^{b}.\]  The graded Betti numbers of $U \cap V$ 
can be computed using Theorem \ref{betti numbers many J}.

We now generalize the proof in \cite{Fra} to compute the graded Betti numbers 
of $U$ and $V$.

\begin{theorem} \label{betti numbers U}
With the notation as above, for all $i \geq 0$,
{\footnotesize\[
\beta_{i,i+a}(U) = \binom{|C|+a -1}{a + i }\binom{a + i -1}{i} + 
\sum_{j=1}^{a-b}\sum_{k=0}^{|A|-1} \binom{k+j-1}{j-1}\binom{|C|+a-j-1}{a-j}\binom{|C|+k}{i}\]}and 
$\beta_{i,j}(U)=0$ for all other $i,j \geq 0$.
\end{theorem}
\begin{proof}
To compute the graded Betti numbers of $U$, first note that
we know the graded Betti numbers of $\mf_C^a$ by Lemma \ref{betti numbers
one J}.  We shall add the remaining
generators of $U$ to $\mf_C^a$, one at a time, and at each intermediate step,
compute the graded Betti numbers
 of the resulting ideal using Lemma \ref{mapping cone
lemma}.  After adding the last generator, we will arrive at the desired
formula.

We add the remaining generators of $U$ to $\mf_C^a$ in the following
order: First, we add the generators of $\mf_A\mf_C^{a-1}$, then
those of $\mf_A^2\mf_C^{a-2}$, and so on.  When adding the 
generators of $\mf_A^t\mf_C^{a-t}$, we shall add the generators
in descending lexicographic order with respect to the ordering $x_1 > \cdots
> x_{t_1} > z_1 > \cdots > z_{t_3}$.  Let $m_{\ell}$ denote
the $\ell$-th monomial added to $\mf_C^a$, and set $U_{\ell} = \mf_C^a +
(m_1,\ldots,m_{\ell})$. 

For each $m = x_1^{a_1}\cdots x_{t_1}^{a_{t_1}}z_1^{c_1}\ldots z_{t_3}^{c_{t_3}}\in \mf_A^t\mf_C^{a-t}$, 
we associate to $m$ the following number:
\[
k_x(m) := \max\{k ~|~ x_{k+1} ~\mbox{divides}~ x_1^{a_1}\cdots x_{t_1}^{a_{t_1}}\}.\]
For example $k_x(x_1^4z_1^{c_1}\ldots z_{t_3}^{c_{t_3}})=0$ since $x_1$ divides $x_1^4$, while
 $k_x(x_1^2x_2x_3z_1^{c_1}\ldots z_{t_3}^{c_{t_3}}) = 2$
because $x_3$ divides $x_1^2x_2x_3$.  
This notation is needed to prove:
\vspace{.25cm}

\noindent
{\bf Claim:} If $m_{\ell}$, the $\ell$-th monomial to be added $\mf_C^a$,
belongs to $m_{\ell} \in \mf_A^t\mf_C^{a-t}$ and $k=k_x(m_{\ell})$ then
\[U_{\ell-1}:m_{\ell} = \mf_C + (x_1,\ldots,x_{k}).\]

\begin{proof}
By construction, 
\[U_{\ell-1} = \mf_C^a + \mf_A\mf_C^{a-1} + \cdots + \mf_A^{t-1}\mf_C^{a-t+1} 
+ ( m \in \G(\mf_A^t\mf_C^{a-t}) ~|~ m > m_{\ell} ).\]
Since multiplying $m_{\ell}$ by any $z_i$ gives $z_im_{\ell}
\in \mf_A^t\mf_C^{a-t+1} \subseteq \mf_A^{t-1}\mf_C^{a-t+1}$, it immediately
follows that $\mf_C \subseteq U_{\ell-1}:m_{\ell}$.  

If $k =0$, then $m_{\ell} = x_1^tz_1^{c_1}\cdots z_{t_3}^{c_{t_3}}$.  Multiplying $m_{\ell}$ by
any monomial $m \in R$ not divisible by $z_i$ does not land you in 
$\mf_C^{a-i}$ for $i=0,\ldots,t-1$.  So, if $mm_{\ell} \in U_{\ell-1}$, then $mm_{\ell}
\in \mf_A^t\mf_C^{a-t}$.  That is $mm_{\ell}$ must be divisible by
a monomial in $\mf_A^t\mf_C^{a-t}$ greater than $m_{\ell}$.  But the only elements greater
than $m_{\ell}$ must have the form  $x_1^tz_1^{d_1}\cdots z_{t_3}^{d_{t_3}}$ with
$z_1^{d_1}\cdots z_{t_3}^{d_{t_3}} > z_1^{c_1}\cdots z_{t_3}^{c_{t_3}}$.  No element
of this form can divide $mm_{\ell}$.
So, if $k=0$, $U_{\ell-1}:m_{\ell} = \mf_C$.

If $k > 0$, to show
that $x_1,\ldots,x_k \in U_{\ell-1}:m_{\ell}$, we note
that $m_{\ell} = x_1^{a_1}\cdots x_{k+1}^{a_{k+1}}z_1^{c_1}\cdots z_{t_3}^{c_{t_3}}$. 
Then for each $i = 1,\ldots,k$, 
\[x_i m_{\ell} = (x_ix_1^{a_1} \cdots x_{k+1}^{a_{k+1}-1}
z_1^{c_1}\cdots z_{t_3}^{c_{t_3}})x_{k+1} = m'x_{k+1}.\]
Now $m' > m_{\ell}$ with respect to our ordering, so 
$x_im_{\ell} \in U_{\ell-1}$.  So $\mf_C + (x_1,\ldots,x_k) \subseteq
U_{\ell-1}:m_{\ell}$. 

To prove the reverse inclusion, let $m$ be any monomial of $R$
not divisible by either the $z_i$s or $x_1,\ldots,x_k$.  If
$mm_{\ell} \in U_{\ell-1}$, then $mm_{\ell}$ must be in
$\mf_A^t\mf_C^{a-t}$ since $mm_{\ell} \not\in \mf_A^i\mf_C^{a-i}$
for $i=0,\ldots,t-1$.  For $mm_{\ell}$ to be both in $U_{\ell-1}$
and $\mf_A^t\mf_C^{a-t}$, it must be divisible by some monomial $m' \in \mf_A^t\mf_C^{a-t}$
with  $m' > m_{\ell}$. For $m'$ to be larger than $m_{\ell} = x_1^{a_1}\cdots x_{k+1}^{a_{k+1}}z_1^{c_1}
\cdots z_{t_3}^{c_{t_3}}$, the exponent of one of $x_1,\ldots,x_{k+1},z_1,\ldots,z_{t_3}$ must be larger in $m'$. Let $g$ be the first index where the exponent of some $x$ or $z$ variable is bigger in $m'$ than in $m_{\ell}$.

We claim that $g \not = k+1$. Since $m' \in \G(\mf_A^t\mf_C^{a-t})$, we can write $m'$ as \[ m'= x_1^{b_1} \cdots x_{t_1}^{b_{t_1}} z_1^{d_1} \cdots z_{t_3}^{d_{t_3}},\] where $b_1+\cdots+b_{t_1}=a$. If $g=k+1$, then $b_{k+1} > a_{k+1}$. By the definition of $g$, $a_i = b_i$ for $i=1,\dots,k$. But then we have \[ a = a_1 + \cdots + a_{k+1} < b_1 + \cdots + b_{k+1} \le b_1 + \cdots + b_{t_1} = a,\] a contradiction.

Hence $m'$ is divisible either by some $z_i$ or one of $x_1,\ldots,x_k$, and thus, $m$ would also have this property, providing us with a contradiction. So the only monomials in $U_{\ell-1}:m_{\ell}$ are those in $\mf_C + (x_1,\ldots,x_k)$.
\end{proof}

We now compute the graded Betti numbers of $U_{\ell}$ for each $\ell$.
When $\ell = 0$, $U_{0} = \mf_C^a$, and 
the graded Betti numbers are given by Lemma \ref{betti numbers one J}:
\[\beta_{i,i+a}(U_0) = \binom{|C|+a-1}{a+i}\binom{a+i-1}{i}\]
and $\beta_{i,j}(U_0)=0$ for all other $i,j \geq 0$.  Observe that this
formula implies that the regularity of $R/U_0$ is $a-1$.

Suppose now that $\ell > 0$, and that $m_{\ell}$ is the $\ell$-th monomial.
Furthermore, suppose that $m_{\ell} \in \mf_A^t\mf_C^{a-t}$ with
$k=k_x(m_{\ell})$.  Applying the claim, we have a short exact sequence
\[0 \rightarrow R/(U_{\ell-1}:m_{\ell})(-a) =
R/(\mf_C + (x_1,\ldots,x_k))(-a) \stackrel{\times m_{\ell}}{\longrightarrow}
 R/U_{\ell-1} \rightarrow R/U_{\ell} \rightarrow 0.\]
By Lemma \ref{mapping cone lemma}, the mapping construction gives
a minimal graded resolution of $R/U_{\ell}$. 
Thus
\[\beta_{i,i+a}(U_{\ell}) = \beta_{i,i+a}(U_{\ell-1}) + \binom{|C|+k}{i}.\]
and $\beta_{i,j}(U_{\ell}) = 0$ for all other $i,j \geq 0$.
So, each new generator $m$ that we add to
$U_{0}$ contributes $\binom{|C|+k_x(m)}{i}$ to $\beta_{i,i+a}(U)$.  

For each $t = 1,\ldots,a-b$,  there are $\binom{|C|+a-t-1}{a-t}$
generators of $\mf_A^t\mf_C^{a-t}$ with $k_x(m) = 0$.  These
are the elements of $x_1^t\mf_C^{a-t}$.    Also,
for each $t = 1,\ldots,a-b$, there are 
\[\binom{k+t-1}{t-1}\binom{|C|+a-t-1}{a-t}\] generators
of $\mf_A^t\mf_C^{a-t}$ with $k_x(m) =k$ as $1 \leq k \leq |A|-1$.  To
see this, we first need to count the number of elements of $\mf_A^t$ 
of the form $x_1^{a_1}\cdots x_{k+1}^{a_{k+1}}$ with $a_{k+1} \geq 1$.
This is equivalent to counting the number of nonnegative integer solutions to
\[a_1 + \cdots + a_{k+1} = t ~~\mbox{with $a_{k+1} > 0$}.\]
Standard techniques in combinatorics imply that this equals $\binom{k+t-1}{t-1}$.
For each monomial $m \in \mf_A^t$ of this form, every monomial 
$m'' \in m\mf_C^{a-t}$ has $k_x(m'')=k$.  So we get 
$\binom{k+t-1}{t-1}\binom{|C|+a-t-1}{a-t}$ generators with $k_x(m) =k$.
By the discussion in the previous paragraph, each generator contributes $\binom{|C|+k}{i}$
to $\beta_{i,i+a}(U)$.
The formula in the statement of the theorem then
comes by summing over all $t$ and $k$.
\end{proof}

Note that when $K \subseteq J$ as in Case 3, $C = K \cap J = K$ and $A = J \backslash K$.  So
$I = U$ when $K \subseteq J$.  The above theorem provides the following formula for $I$ in Case 3.

\begin{corollary}
Suppose $J,K \subseteq [n]$ with $K \subseteq J$.  If $I = \mf_{J}^a\cap\mf_K^b$, then
{\footnotesize
\begin{eqnarray*}\beta_{i,i+a}(I) &=& \binom{|K|+a -1}{a + i }\binom{a + i -1}{i} + 
\sum_{j=1}^{a-b}\sum_{k=0}^{|J \backslash K|-1} \binom{k+j-1}{j-1}\binom{|K|+a-j-1}{a-j}\binom{|K|+k}{i}
\end{eqnarray*}}
and $\beta_{i,j}(I) = 0$ for all other $i,j \geq 0$.
\end{corollary}

The formula for the graded Betti numbers of $V$ is proved similarly.

\begin{theorem} \label{betti numbers V}
With the notation as above, for $i \geq 0$ and $2 \leq j \leq b$,
{\footnotesize
\begin{eqnarray*}
\beta_{i,i+a+1}(V) & = &
\sum_{i_1+i_2+i_3 =i} \binom{|B|}{1+i_1}
\binom{|A|+a-b}{a-b+1+i_2}\binom{a-b+i_2}{i_2}\binom{|C|+ b -2}{b-1+i_3}
\binom{b+i_3 -2}{i_3} \\
\beta_{i,i+a+j} (V) & = &
\sum_{k_1=0}^{|A|-1} \sum_{k_2=0}^{|B|-1}\binom{|C|+b-j-1}{b-j}\binom{k_1+a-b+j-1}{a-b+j-1}
\binom{k_2+j-1}{j-1}\binom{|C|+k_1+k_2}{i}
\end{eqnarray*}}
\end{theorem}

\begin{proof}
Set $V_0 = \mf_B\mf_A^{a-b+1}\mf_C^{b-1}$.  We add the remaining generators of 
$V$ to $V_0$, one at a time,
and after adding a new generator, we compute the graded Betti numbers
of the resulting ideal.

We shall add the remaining generators of $V$ in the following order:
First, we add the generators of $\mf_B^2\mf_A^{a-b+2}\mf_C^{b-2}$,
then those of $\mf_B^3\mf_A^{a-b+3}\mf_C^{b-3}$, and so on.  When
adding the generators of $\mf_B^t\mf_A^{a-b+t}\mf_C^{b-t}$, we 
will add them in lexicographic order with respect to 
$y_1 > \cdots > y_{t_2} > x_1 > \cdots >x_{t_1} > z_1 > \cdots z_{t_3}$.
We let $m_{\ell}$ denote the $\ell$-th monomial added to $V_0$, and 
define $V_{\ell} := V_0 + (m_1,\ldots,m_{\ell})$.  

To each monomial $m = y_1^{b_1}\cdots y_{t_2}^{b_{t_2}}x_1^{a_1}\cdots
x_{t_1}^{a_{t_1}}z_1^{c_1}\cdots z_{t_3}^{c_{t_3}} \in \mf_B^t\mf_A^{a-b+t}\mf_C^{b-t}$
we associate the following two numbers:
\begin{eqnarray*}
k_x(m) &:=& \max\{k ~|~ x_{k+1} ~\mbox{divides}~ x_1^{a_1}\cdots x_{t_1}^{a_{t_1}}\}. \\
k_y(m) &:=& \max\{k ~|~ y_{k+1} ~\mbox{divides}~ y_1^{b_1}\cdots y_{t_2}^{b_{t_2}}\}.
\end{eqnarray*}
Using this notation, we shall prove:
\vspace{.25cm}

\noindent
{\bf Claim:} Suppose that $m_{\ell}$ is the ${\ell}$-th monomial
added to $V_0$, and that $m_{\ell} \in \mf_B^t\mf_A^{a-b+t}\mf_C^{b-t}$
with $k_y = k_y(m_{\ell})$ and $k_x = k_x(m_{\ell})$.  Then
\[ V_{\ell-1}:m_{\ell} = \mf_C + (x_1,\ldots,x_{k_x},y_1,\ldots,y_{k_y}).\]

\begin{proof}  By definition
\[V_{\ell-1} = \mf_B\mf_A^{a-b+1}\mf_C^{b-1}
+ \cdots +\mf_B^{t-1}\mf_A^{a-b+t-1}\mf_C^{b-t+1} +
( m \in \G(\mf_B^t\mf_A^{a-b+t}\mf_C^{b-t}) ~|~ m > m_{\ell}).\]
It is straightforward to check that $\mf_C \subseteq V_{\ell-1}:m_{\ell}$.

If $k_x = k_y = 0$, then $m_{\ell} = y_1^tx_1^{a-b+t}m'$ where
$m' \in \mf_C^{b-t}$. 
Multiplying $m_{\ell}$ by
any monomial $m \in R$ not divisible by $z_i$ does not land you in 
$\mf_C^{b-i}$ for $i=1,\ldots,t-1$.  So, if $mm_{\ell} \in U_{\ell-1}$, then $mm_{\ell}
\in \mf_B^t\mf_A^{a-b+t}\mf_C^{b-t}$.  That is, $mm_{\ell}$ must be divisible by
a monomial in $\mf_B^t\mf_A^{a-b+t}\mf_C^{b-t}$ greater than $m_{\ell}$.  But the only elements greater
than $m_{\ell}$ must have the form  $y_1^tx_1^{a-b+t}m''$ with
$m'' > m'$.  No element
of this form can divide $mm_{\ell}$.
So, if $k_x = k_y =0$, $V_{\ell-1}:m_{\ell} = \mf_C$.

If $k_y > 0$, then $m_{\ell} = y_{1}^{b_1}\cdots y_{k_y+1}^{b_{k_y+1}}
x_{1}^{a_{1}}\cdots x_{t_1}^{a_{t_1}}z_1^{c_1}\cdots z_{t_3}^{c_{t_3}}$.
Then for each $i = 1,\ldots,k_y$,
\[y_im_{\ell} = (y_iy_{1}^{b_1}\cdots y_{k_y+1}^{b_{k_y+1}-1}
x_{1}^{a_{1}}\cdots x_{t_1}^{a_{t_1}}z_1^{c_1}\cdots z_{t_3}^{c_{t_3}})y_{k_y+1}
= m'y_{k_y+1}\]
But $m' > m_{\ell}$, so $m' \in V_{\ell-1}$, thus 
$y_i \in V_{\ell-1}:m_{\ell}$.  If $k_x > 0$, a similar argument
implies that $x_1,\ldots,x_{k_x} \in V_{\ell-1}:m_{\ell}$.
Hence $\mf_C + (x_1,\ldots,x_{k_x},y_1,\ldots,y_{k_y}) \subseteq
V_{\ell-1}:m_{\ell}$.

The opposite containment follows from an argument similar to the one in Theorem~\ref{betti numbers U}.
\end{proof}

We now compute the graded Betti numbers of $V_{\ell}$ for each $\ell$.
When $\ell = 0$, $V_0 = \mf_B\mf_A^{a-b+1}\mf_C^{b-1}$.  The graded
Betti numbers follow from Theorem \ref{betti numbers many J}:
{\footnotesize
\[\beta_{i,i+a+1}(V_0) = \sum_{i_1+i_2+i_3 =i} \binom{|B|}{1+i_1}
\binom{|A|+a-b}{a-b+1+i_2}\binom{a-b+i_2}{i_2}\binom{|C|+ b -2}{b-1+i_3}
\binom{b+i_3 -2}{i_3} \]}
and $\beta_{i,j}(V_0) = 0$ otherwise.

Suppose that $\ell > 0$ and let $m_{\ell}$ be the $\ell$-th monomial
with $m_{\ell} \in \mf_B^t\mf_A^{a-b+t}\mf_C^{b-t}$.  
We have the short exact sequence
\[0 \rightarrow R/(V_{\ell-1}:m_{\ell})(-a-t)
\stackrel{\times m_{\ell}}{\longrightarrow} R/V_{\ell-1} 
\rightarrow R/V_{\ell} \rightarrow 0.\]
Note that reg$(R/V_0)=a$, and inductively, for $\ell-1 \geq0$, reg$(R/V_{\ell-1}) = a+t-1$ since $V_{\ell-1}:m_{\ell}$ is generated by a subset of the variables. Therefore by Lemma \ref{mapping cone lemma}, the mapping cone construction
gives a minimal graded free resolution of $R/V_{\ell}$. If 
$k_x = k_x(m_{\ell})$ and $k_y = k_y(m_{\ell})$, then the
claim implies $R/(V_{\ell-1}:m_{\ell}) = R/(\mf_C + (x_1,\ldots,x_{k_x},
y_1,\ldots,y_{k_y}))$.  So, each generator $ m \in 
\mf_B^t\mf_A^{a-b+t}\mf_C^{b-t}$
contributes $\binom{|C|+k_x(m) + k_y(m)}{i}$ to $\beta_{i,i+a+t}(V)$.

Counting as in Theorem \ref{betti numbers U} and summing over all possible $t$, $k_x$, and $k_y$, we obtain the final formulas; we leave the details to the reader.
\end{proof}

\begin{theorem}\label{resolution theorem}
Suppose $J,K \subseteq [n]$ are such that $J \cap K, J \backslash K, K \backslash J \neq \emptyset$.
If $I = \mf_J^{a} \cap \mf_K^{b}$, then
\begin{eqnarray*}
\beta_{i,i+a}(I) &=& \beta_{i,i+a}(U) \\
\beta_{i,i+a+1}(I) & = & \beta_{i,i+a+1}(V) + \beta_{i-1,i+a+1}(\mf_B\mf_A^{a-b+1}\mf_C^b)\\
\beta_{i,i+a+j}(I) & = & \beta_{i,i+a+j}(V) ~~\mbox{for $j=2,\ldots,b$}.
\end{eqnarray*}
where $U$ and $V$ are as defined above.
\end{theorem}
\begin{proof}
Since $I = U+V$ is a splitting, the formulas are a consequence of Theorem \ref{EKF} and the
fact that $\beta_{i-1,i+a+1}(U \cap V) = \beta_{i-1,i+a+1}(\mf_B\mf_A^{a-b+1}\mf_C^b)$.
\end{proof}

\section{Applications: multiplicity, combinatorics, and fat points in multiprojective space} \label{s.fatpoints}

In our final section, we present some applications of our results in the earlier sections. First, we discuss some 
cases of the Multiplicity Conjecture of Herzog, Huneke, and Srinivasan. In addition, we use our componentwise 
linearity results and Alexander duality to prove a corollary about the sequential Cohen-Macaulayness of some simplicial 
complexes. Finally, we apply our earlier results to investigate the resolutions of some sets of fat points in 
multiprojective space. The main result of \cite{Fra} is that ideals of small sets of general fat points in 
$\pr^n$ are componentwise linear; we generalize this theorem to multiprojective 
space. Furthermore, we extend work from \cite{Fa,Va,FL,Fra} to describe the graded Betti numbers of ideals of 
small sets of fat points in linear general position in multiprojective space.

\subsection{Multiplicity Conjecture} The Multiplicity Conjecture of Herzog, Huneke, and Srinivasan 
(see, e.g., \cite{HS}) proposes bounds for the multiplicity of an ideal in terms of the shifts in its graded free 
resolution.  The explicit statement is given below.

\begin{conjecture}
Let $R/I$ be a homogeneous $k$-algebra with resolution of the form
\[0 \longrightarrow   \bigoplus_{j=1}^{b_r}
R(-d_{rj}) \longrightarrow \cdots \longrightarrow
\bigoplus_{j=1}^{b_1}
R(-d_{1j}) \longrightarrow R \longrightarrow R/I
\longrightarrow 0.\]
Set $m_i = \min\{ d_{ij} ~|~ j = 1,\ldots,b_i\}$ and $M_i
= \max\{d_{ij} ~|~ j = 1,\ldots,b_i\}$.  If codim$(I) = c$
and $e(R/I)$ denotes the multiplicity of $R/I$, then 
\[ e(R/I) \leq \frac{\prod_{i=1}^c M_i}{c!}.\]
Furthermore, if $R/I$ is Cohen-Macaulay, then
\[\frac{\prod_{i=1}^c m_i}{c!} \leq e(R/I). \]
\end{conjecture}

In \cite{Roemer}, R\"omer proved that when the characteristic of $k$ is zero, componentwise linear 
ideals satisfy the above Multiplicity Conjecture.
As a consequence of Theorem~\ref{t.maintool}, Corollary~\ref{c.2CWL}, Theorem~\ref{t.threeideals}, and 
R\"{o}mer's result, we have:

\begin{corollary}
\label{c:mc}
Suppose $\operatorname{char}(k)=0$. Let $I=\mf_{J_1}^{a_1} \cap \cdots \cap \mf_{J_s}^{a_s}$, and suppose either that $s \le 3$, or $J_i \cup J_j = [n]$ for all $i \not = j$. Then $I$ satisfies the upper bound of the Multiplicity Conjecture.
\end{corollary}

Note that we only know that the upper bound is true since in general, $R/I$ 
may not be Cohen-Macaulay. If it is, then the lower bound holds as well. (R\"{o}mer states his result only for 
the upper bound, but his proof is based on the fact that if $I$ is componentwise linear, then $I$ and the
reverse-lex generic initial ideal $\gin(I)$ have the same graded Betti numbers in characteristic zero. Both 
bounds of the conjecture hold for all Cohen-Macaulay generic initial ideals in characteristic zero since the 
bounds are true for all Cohen-Macaulay strongly stable ideals. Since the reverse-lex gin preserves depth and 
dimension, if $R/I$ is Cohen-Macaulay, $R/\gin(I)$ is as well, so the lower bound holds in that case.)

\subsection{The sequentially Cohen-Macaulay property}
The notion of componentwise linearity is intimately related to the concept of sequential Cohen-Macaulayness.

\begin{definition} \label{d.seqcm}
Let $R=k[x_1,\dots,x_n]$. A graded $R$-module $M$ is called {\bf sequentially Cohen-Macaulay} if there exists a finite filtration of $M$ by graded $R$-modules
\[ 0 = M_0 \subset M_1 \subset \cdots M_r = M \]
such that each $M_i/M_{i-1}$ is Cohen-Macaulay, and the Krull dimensions of the quotients are increasing:
\[\dim (M_1/M_0) < \dim (M_2/M_1) < \cdots < \dim (M_r/M_{r-1}).\]

 We say that a simplicial complex $\Delta$ is sequentially Cohen-Macaulay if $R/I_{\Delta}$ is
sequentially Cohen-Macaulay, where $I_\Delta$ is the Stanley-Reisner ideal of $\Delta$. 
\end{definition}

Stanley introduced sequential Cohen-Macaulayness in connection with developments in the theory of shellability; 
see, e.g., \cite{Stanley} for a definition of shellable. A shellable pure simplicial complex (that is, a shellable 
simplicial complex whose maximal faces all have the same dimension) is Cohen-Macaulay, but if one extends the 
definition of shellability to allow nonpure simplicial complexes, one obtains simplicial complexes that are not 
Cohen-Macaulay. However, they are sequentially Cohen-Macaulay. 

The theorem connecting sequentially Cohen-Macaulayness to componentwise linearity is based on the idea of Alexander duality. We define Alexander duality for squarefree monomial ideals and then state the fundamental result of
 Herzog and Hibi \cite{HH} and Herzog, Reiner, and Welker \cite{HRW}.

\begin{definition} If $I = (x_{1,1}x_{1,2}\cdots x_{1,t_1}, \ldots, x_{s,1}x_{s,2}\cdots x_{s,t_s})$ is a
squarefree monomial ideal, then the {\bf Alexander dual} of $I$, denoted $I^{\star}$, is the monomial
ideal
\[I^{\star} = (x_{1,1},\ldots,x_{1,t_1}) \cap \cdots \cap (x_{s,1},\ldots,x_{s,t_s}).\]
If $\Delta$ is a simplicial complex and $I = I_{\Delta}$ its Stanley-Reisner ideal, then the
simplicial complex $\Delta^{\star}$ with $I_{\Delta^{\star}} = I^{\star}_{\Delta}$ is the {\bf Alexander dual} of $\Delta$.
\end{definition}

\begin{theorem} \label{t.seqcm}
Let $\Delta$ be a simplicial complex with Stanley-Reisner ideal $I_{\Delta}$. Let $\Delta^*$ be the Alexander dual of $\Delta$. Then $R/I_{\Delta}$ is sequentially Cohen-Macaulay if and only if $I^{\star}_{\Delta} = 
I_{\Delta^{\star}}$ is componentwise linear.
\end{theorem}

Our results in this paper yield the following corollary.

\begin{corollary} \label{c.combinatorics}
Let $\Delta$ be a simplicial complex on $n$ vertices, and let $I_{\Delta}$ be its Stanley-Reisner ideal, minimally generated by squarefree monomials $m_1,\dots,m_s$. If $s \le 3$, so that $\Delta$ has at most three minimal nonfaces, or if $\Supp(m_i) \cup \Supp(m_j) = \{x_1,\dots,x_n\}$ for all $i \not = j$, then $\Delta$ is sequentially Cohen-Macaulay.
\end{corollary}

\begin{proof}
$I_\Delta$ is a squarefree monomial ideal; suppose it is minimally generated by monomials 
$\{x_{1,1} \cdots x_{1,t_1},\dots, x_{s,1} \cdots x_{s,t_s}\}$. Then 
\[ I^{\star}_{\Delta} = I_{\Delta^{\star}} = (x_{1,1},\dots,x_{1,t_1}) \cap \cdots \cap (x_{s,1}, \dots, x_{s,t_s}).\]
By Theorem~\ref{t.maintool}, Corollary~\ref{c.2CWL}, or Theorem~\ref{t.threeideals}, $I^{\star}_{\Delta}$ is 
componentwise linear, and so Theorem~\ref{t.seqcm} gives the result.
\end{proof}

\begin{example} \label{e.combinatorics}
Let $\Delta$ be a simplicial complex on six vertices. Suppose the minimal nonfaces of 
$\Delta$ are $\{145,126,135\}$. Then 
\[I_{\Delta^{\star}}=(x_1,x_4,x_5) \cap (x_1,x_2,x_6) \cap (x_1,x_3,x_5) \subset R=k[x_1,\dots,x_6]\]
is componentwise linear by Theorem~\ref{t.threeideals}, and thus $\Delta$ is sequentially Cohen-Macaulay. Note that $\Delta$ is not Cohen-Macaulay since codim $I_{\Delta}=1$, while the projective dimension of $R/I_{\Delta}$ is two.
\end{example}

\subsection{Fat points in multiprojective space} 
We begin by recalling some of the relevant definitions for points in multiprojective space 
(for more on this topic see \cite{VT1,VT2,VT3}).  The coordinate ring of  
$\pr^{n_1} \times \cdots \times \pr^{n_r}$ is
the $\N^r$-graded polynomial ring 
$R = k[x_{1,0},\ldots,x_{1,n_1},\ldots,x_{r,0},\ldots,x_{r,n_r}]$
with $\deg x_{i,j} = e_i$, the $i$-th basis vector of $\N^r$.  The
defining ideal of a point  
$P = P_1 \times \cdots \times P_r \in \pr^{n_1} \times \cdots \times \pr^{n_r}$
is the prime ideal $I_{P} = (L_{1,1},\ldots,L_{1,n_1},\ldots,L_{r,1},\ldots,L_{r,n_r})$ with
$\deg L_{i,j} = e_i$.   The forms $L_{i,1},\ldots,L_{i,n_i}$  are the generators
of the defining ideal of $P_i \in \pr^{n_i}$.

\begin{definition}
A set of points $X \subseteq \pr^n$ is said to be in {\bf linear general position} if no more than 
two points lie on a line, 
no more than three points line in a plane, ..., no more than $n$ points lie in an $(n-1)$-plane.
\end{definition}
Observe that the above definition is equivalent to the fact that if $\LL_d$ is 
any linear subspace of $\pr^n$ of dimension $d$ with $d = 0,\ldots,n-1$, 
then the intersection of $\LL_d$ and $X$ contains at most $d+1$ points of $X$. When $d=0$, $\LL_d$
is a point, so this simply says that the intersection of a point and $X$ is at most one point.
To extend this to a multigraded context, we say that $\LL$ is $(d_1,\ldots,d_r)$-linear
subspace of $\pr^{n_1} \times \cdots \times \pr^{n_k}$ if $\LL = \LL_{d_1} \times 
\cdots \times \LL_{d_r}$, where each $\LL_{d_i}$ is a linear subspace of $\pr^{n_i}$
of dimension $d_i$ with $d_i = 0,\ldots,n_i$ (so  $\LL_{n_i} = \pr^{n_i}$ is allowed) and
there exists at least one $j \in [r]$ such that $d_j < n_j$.  

\begin{definition}
A set of points $X \subseteq \pr^{n_1} \times \cdots 
\times \pr^{n_r}$ is in {\bf linear general position} 
if for every $(d_1,\ldots,d_r)$-linear subspace $\LL$, the intersection
of $\LL$ and $X$ contains at most $d+1$ points of $X$ where $d = \min\{d_1,\ldots,d_r\}$. 
\end{definition}

We point out that if $\LL$ is $(d_1,\ldots,d_r)$-linear subspace with $d= d_i = 0$, then $\LL_{d_i}$ is a point. 
So if $X$ is in linear general position, this means that at most one point of $X$ can
intersect $\LL$, which, in turn, implies that at most one point of $X$ can have $i$-th coordinate equal to $\LL_{d_i}$.  
It follows from this observation that for any two points $P,Q \in X$ with $X$ in linear
general position in $\pr^{n_1} \cdots \times \cdots \pr^{n_k}$, we must have
$P_i \neq Q_i$ for $i = 1,\ldots,r$.  
In other words if $\pi_i:\pr^{n_1} \times \cdots \times \pr^{n_r} \rightarrow \pr^{n_i}$ denotes the 
projection morphism for $i = 1,\ldots,r$, and if $\{Q_1,\dots,Q_t\} \in \pr^{n_1} \times \cdots \times \pr^{n_r}$ 
is in  linear general position, then the sets of the projections $\{\pi_i(Q_1), \dots, \pi_i(Q_t)\}$ are in 
linear general position in $\pr^{n_i}$ for each $i$.   In particular, we require that 
$\pi_i(Q_j) \not = \pi_i(Q_l)$ for all $i$ and all $j \not = l$; 
see Remark~\ref{r:notgeneral} for what can go wrong without this condition. 

\begin{definition} Let $\{P_1,\ldots,P_s\} \subseteq \pr^{n_1} \times \cdots \times \pr^{n_r}$ be
a set of points with the defining ideal of $P_i$ denoted $I_{P_i}$ and let $a_1,\ldots,a_s$ be positive integers.
The scheme $Z \subseteq \pr^{n_1} \times \cdots \times \pr^{n_r}$ defined by
\[I_Z = I_{P_1}^{a_1} \cap \cdots \cap I_{P_s}^{a_s}\]
is scheme of {\bf fat points}, and is sometimes denoted $Z = \{(P_1,a_1),\ldots,(P_s,a_s)\}$.
We call $a_i$ the {\bf multiplicity} of the point $P_i$.  The points $\{P_1,\ldots,P_s\}$
are referred to as the {\bf support} of $Z$.
\end{definition}

By a small set of linear general fat points in $\pr^n$, 
we mean that the support has at most $n+1$ points in linear general position.
This restriction allows us to make a change of coordinates to move all the points to the coordinate 
vertices, and we can take the ideal corresponding to the set of fat points to be an intersection of 
monomial ideals 
\[ I=(x_1,\dots,x_n)^{a_0} \cap (x_0,x_2,\dots,x_n)^{a_1} \cap \cdots \cap (x_0,\dots,x_{s-1},x_{s+1},x_n)^{a_s} .\]

If we are working in $\pr^{n_1} \times \cdots \times \pr^{n_r}$,
we would like to change coordinates to work with a set of fat points 
at the coordinate vertices so that the corresponding ideals are monomial ideals. Therefore, a small 
set of fat points can consist of no more than $1+\min \{n_1, \dots, n_r \}$ points.  The set
of fat points is general if the points in the support are in linear general position.

Suppose that $I$ is the ideal of a small set of general fat points in 
$\pr^{n_1} \times \cdots \times \pr^{n_r}$. As a consequence of Theorem~\ref{t.maintool}, we can generalize 
the componentwise linearity result for the $r=1$ case from \cite{Fra} 
(and obtain a different proof for that case).

\begin{theorem}\label{t:smallfatcwl}
Let $I$ be the ideal of $s+1$ general fat points in $\pr^{n_1} \times \cdots \times \pr^{n_r}$, where $s \le \min \{n_1, \dots, n_r \}$. Then for all $d$, $(I_d)$ is polymatroidal, and $I$ is componentwise linear.
\end{theorem}

\begin{proof}
Because $I$ is the ideal of a small set of general fat points in multiprojective space, we may assume that 
$I$ has the form 
\[ I=(x_{1,1},\dots,x_{1,n_1},x_{2,1},\dots,x_{2,n_2},\dots,x_{r,1},\dots,x_{r,n_r})^{a_0} \cap \cdots \cap \] 
\[ (x_{1,0},\dots,\hat{x_{1,s}},\dots, x_{1,n_1},x_{2,0},\dots, \hat{x_{2,s}}, \dots, x_{2,n_2},\dots,x_{r,0},\dots, 
\hat{x_{r,s}}, \dots, x_{r,n_r})^{a_s} \subset R, \] 
where $\hat{x_{i,s}}$ denotes that $x_{i,s}$ is left out. Note that the union of the variables appearing in any 
two of the components is all the variables of $R$. Hence the result follows immediately from Theorem~\ref{t.maintool}.
\end{proof}

As in Corollary~\ref{c:mc}, when the char$(k)=0$, 
Theorem~\ref{t:smallfatcwl} implies that ideals of small sets of general fat points in multiprojective space 
satisfy the Multiplicity Conjecture of Herzog, Huneke, and Srinivasan. Note that if $r > 1$, the ideal will 
not be Cohen-Macaulay
(for example, see \cite{VT1,VT2}), so we may only conclude that the conjectured upper bound is true. 

We conclude this discussion with a remark about how we defined the notion of a ``general'' set of fat points.

\begin{remark}\label{r:notgeneral}
In our definition of what it means for a set of fat points $Q_1, \dots, Q_s$ in multiprojective space to be general, 
we required that 
for all $i$ and all $j \not = l$, the projections $\pi_i(Q_j) \not = \pi_i(Q_l)$. If that condition 
is not satisfied, the corresponding ideal may not be componentwise linear.

Consider the points $[1:0]\times [1:0],[1:0] \times [0:1], [0:1]\times [1:0],$ and  
$[0:1]\times [0:1]$ in 
$\pr^1 \times \pr^1$, and suppose each point has multiplicity one. The ideal corresponding 
to the set of four points in $R = k[x_0,x_1,y_0,y_1]$ is 
\[ I = (x_1,y_1) \cap (x_1,y_0) \cap (x_0,y_1) \cap (x_0,y_0) = (x_0x_1, y_0y_1). \] 
This ideal is a complete intersection of degree two polynomials, and hence it is not componentwise 
linear; in particular, $I=(I_2)$ does not have a linear resolution. The problem is that the union of 
the variables appearing in, for example, the first two components, is not all of $\{x_0,x_1,y_0,y_1\}$.
\end{remark}

We turn now to the graded Betti numbers of two general fat points in multiprojective space. 
As an application of Theorem \ref{resolution theorem} we can compute the
$\N$-graded Betti numbers of the defining ideal of two fat points in 
$\pr^{n_1} \times \cdots \times \pr^{n_r}$ in linear general position.  

\begin{corollary} \label{resolution two fat points}
Let $Z = \{(P,a),(Q,b)\}$ be two fat points in $\pr^{n_1} \times \cdots \times \pr^{n_r}$
with $a \geq b$
Set $N = n_1 + \cdots + n_r$, and let $I_Z$ denote the defining ideal of $Z$.  
If $P$ and $Q$ are in linear general position, then
\footnotesize
\begin{eqnarray*}
\beta_{i,i+a}(I_Z) &=& 
\binom{N-r+a -1}{a + i }\binom{a + i -1}{i} + 
\sum_{j=1}^{a-b}\sum_{k=0}^{r-1} \binom{k+j-1}{j-1}\binom{N-r+a-j-1}{a-j}\binom{N-r+k}{i}\\
\beta_{i,i+a+1}(I_Z) & = & \sum_{i_1+i_2+i_3 =i} \binom{r}{1+i_1}
\binom{r+a-b}{a-b+1+i_2}\binom{a-b+i_2}{i_2}\binom{N-r+ b -2}{b-1+i_3}
\binom{b+i_3 -2}{i_3} +\\
&& \sum_{i_1+i_2+i_3 =i-1} \binom{r}{1+i_1}
\binom{r+a-b}{a-b+1+i_2}\binom{a-b+i_2}{i_2}\binom{N-r+ b -1}{b+i_3}
\binom{b+i_3 -1}{i_3}
\\
\beta_{i,i+a+j}(I_Z) & = & 
\sum_{k_1=0}^{r-1} \sum_{k_2=0}^{r-1}\binom{N-r+b-j-1}{b-j}\binom{k_1+a-b+j-1}{a-b+j-1}
\binom{k_2+j-1}{j-1}\binom{N-r+k_1+k_2}{i}
\\
&& \mbox{for $j=2,\ldots,b$.} 
\end{eqnarray*}
\normalsize
and $\beta_{i,j}(I_Z) = 0$ for all other $i,j \geq 0$.
\end{corollary}

\begin{proof}
Since $P$ and $Q$ are in linear general position, we may assume (after a change of
coordinates) that $P = [1:0:\cdots:0] \times
\cdots \times [1:0:\cdots:0]$ and $Q = [0:1:0:\cdots:0] \times \cdots
\times [0:1:0\cdots:0]$.  So, the defining ideal of $I_Z$ has the form
\[I_Z = (x_{1,1},\ldots,x_{1,n_1},\ldots,x_{r,1},\ldots,x_{r,n_r})^a \cap
(x_{1,0},x_{1,2},\ldots,x_{1,n_1},\ldots,x_{r,0},x_{r,2},\ldots,x_{r,n_r})^b\]
The graded Betti numbers of $I_Z$ are then a consequence of Theorem \ref{resolution theorem} with
$|C| = N -r$ and $|A| = |B| = r$.
\end{proof}

\begin{remark} When $r=1$ in the previous corollary, we recover the formulas of Valla \cite{Va}
and first author \cite{Fra} for two fat points in $\pr^n$.  
When $r > 1$, then $I_Z$ also has a multigraded resolution of the form
\[0 \rightarrow \bigoplus_{\underline{j} \in \N^r} R(-\underline{j})^{\beta_{h,\underline{j}}(I_Z)}
\rightarrow   \cdots \rightarrow \bigoplus_{\underline{j} \in \N^r} R(-\underline{j})^{\beta_{0,\underline{j}}(I_Z)}
\rightarrow I_Z \rightarrow 0.\]
Corollary \ref{resolution two fat points} gives us some information on the 
multigraded Betti numbers $\beta_{i,\underline{j}}(I_Z)$ 
because of the identity $\beta_{i,j}(I_Z) =  \sum_{|\underline{j}| = j} \beta_{i,\underline{j}}(I_Z)$.
\end{remark}


\end{document}